\documentclass[a4paper,UKenglish,numberwithinsect,cleveref]{lipics-v2021}

\usepackage[utf8]{inputenc}
\usepackage{mathtools,microtype,todonotes}
\hypersetup{hidelinks}

\theoremstyle{definition}
\newtheorem{problem}[theorem]{Problem}
\crefname{problem}{Problem}{Problems}
\Crefname{problem}{Problem}{Problems}
\crefname{claim}{Claim}{Claims}
\Crefname{claim}{Claim}{Claims}
\theoremstyle{plain}

\nolinenumbers

\newcommand{\N}{\mathbb N}
\newcommand{\Z}{\mathbb Z}
\newcommand{\Q}{\mathbb Q}
\newcommand{\W}{\mathsf W}
\newcommand{\poly}{\mathsf{poly}}
\newcommand{\Hall}{\mathsf{Hall}}
\newcommand{\pSSS}{\mathsf{pSSS}}
\newcommand{\WCS}{\mathsf{WCS}}
\newcommand{\pF}{\mathsf{pF}}
\newcommand{\SolPG}[1]{\mathsf{Sol}_{#1}\mathsf{PG}}
\newcommand{\NilPG}[1]{\mathsf{Nil}_{#1}\mathsf{PG}}
\newcommand{\dig}{\mathsf{dig}}
\DeclareMathOperator{\Sym}{Sym}
\newcommand{\wCNFSAT}{\mathsf{wCNFSAT}}
\newcommand{\FPT}{\mathsf{FPT}}
\newcommand{\fpt}{\mathsf{fpt}}
\newcommand{\func}{\mathsf{func}}
\newcommand{\lfpt}{\leq^{\fpt}}
\newcommand{\rest}{\mathord\restriction}

\title{Parameterized Subset Sum for Integers and Solvable Permutation Groups}
\titlerunning{Subset sum for integers and solvable permutation groups}
\author{Markus Lohrey}{Universit\"at Siegen, Germany}{}{https://orcid.org/0000-0002-4680-7198}{}
\authorrunning{M.~Lohrey}
\Copyright{}
\ccsdesc[500]{Theory of computation~Parameterized complexity and exact algorithms}
\keywords{Subset sum, parameterized complexity, permutation groups, solvable groups, nilpotent groups}

\hideLIPIcs

\begin{document}
\maketitle

\begin{abstract}
We show that parameterized subset sum for binary-encoded integers belongs
to $\W[2]$, where the parameter is the number of selected integers. This improves a result of Buss and Islam,
who proved membership in $\W[3]$.
For every fixed integer $d\geq1$, we also show that parameterized subset
sum for permutation groups of derived length $d$ belongs to
$\W[2d]$. The first proof gives an fpt-reduction to weighted CNF
satisfiability without changing the parameter. The second proof combines
positive CNFs for addition in small abelian groups with an affine
decomposition of a solvable permutation action. This gives circuits of
weft at most $2d$ and depth depending only on $d$ for product evaluation.
For every fixed nilpotency class, we obtain the stronger bound $\W[2]$
by an fpt-reduction to subset sum in abelian permutation groups.
The reduction uses Hall's polynomial coordinates in a free nilpotent
group and the Magnus embedding.
For nilpotent permutation groups without a fixed class bound, we prove
$\W[2]$-hardness, already for permutation $2$-groups.
\end{abstract}

\section{Introduction}

The classical subset sum problem asks whether some of the integers in a given
list sum to a given target integer. With binary-encoded integers, it is one of the
fundamental $\mathsf{NP}$-complete problems \cite{Karp1972}.
Its algorithmic complexity continues to attract considerable attention.
For positive integers, Bringmann \cite{Bringmann2017} gave a randomized
pseudopolynomial algorithm whose running time is near-linear in the number
of input integers plus the target integer, up to logarithmic factors.
Chan \cite{Chan2026} recently obtained a deterministic algorithm with the
same bound, up to logarithmic factors.
Independent recent preprints by Chen, Mao, and Zhang
\cite{ChenMaoZhang2026} and Jin and Mao \cite{JinMao2026} give
deterministic pseudopolynomial algorithms with almost-linear dependence
on the largest absolute value of the input integers.

In this paper, we study the parameterized complexity of subset sum, where
the parameter $k\in\N$ is the number of selected integers from the input list.
Parameterized algorithms aim to confine the expensive part of a computation
to a parameter that is determined by the input but may be small compared to the total input length.
A problem belongs to the class $\FPT$ (for fixed-parameter tractable) if it can be solved with an algorithm
that for input length $n\in\N$ has a running time of $f(k)\cdot n^c$, where $f:\N\to\N$ is computable
and $c\in\N$ is a constant independent of $k$. 
Attempts to define an analogue of $\mathsf{NP}$ in the parameterized
world led to the $\W$-hierarchy, $\W[1]\subseteq\W[2]\subseteq\cdots$.
Hardness for one of these classes gives evidence against
fixed-parameter tractability, whereas membership results relate a problem
to standard problems in the $\W$-hierarchy such as parameterized clique or weighted satisfiability.
We refer to the books of Downey and Fellows \cite{DowneyFellows1999}
and Flum and Grohe \cite{FlumG06} for the foundations of this theory.

For parameterized subset sum over the integers, denoted by $\pSSS[\Z]$ in this paper,
the best previously known unconditional bounds are $\W[1]$-hardness,
proved by Downey and Fellows \cite[Theorem~4.3]{DowneyFellows1995II},
and membership in $\W[3]$, proved by Buss and Islam \cite{BussIslam2007}.
Abboud, Lewi, and Williams \cite{AbboudLewiWilliams2014} proved
$\W[1]$-completeness for the zero-target version when the $m\geq1$ input
integers have absolute value at most $m^{2k}$.
For arbitrary integers, they obtained a randomized fpt-reduction to
parameterized clique. They also showed that this reduction can be
derandomized under a circuit lower-bound assumption, saying that there is a Boolean function
computable in time $2^{\mathcal{O}(n)}$ on inputs of length $n$ that requires
circuits of size $2^{\Omega(n)}$ even when oracle gates for SAT are
allowed. Under this assumption, $\pSSS[\Z]$ is $\W[1]$-complete;
see Appendix~B of the full version of \cite{AbboudLewiWilliams2014}.

Our first result improves the unconditional upper bound for $\pSSS[\Z]$ to $\W[2]$
(\cref{thm:integers}).
After shifting the input integers to make them nonnegative and adjusting the target accordingly, we write them in base $k+1$. 
Since every carry is smaller than $k+1$, in a sum of $k$ selected input integers that yields the target integer,
 the carry into a digit position $p$ is uniquely determined by the sum of the position-$p$ digits of the 
 $k$ selected integers and the position-$p$ digit of the target integer (see \cref{lem:integer-carries}). 
 This leads to a locally checkable condition that can be described by a set 
 of forbidden digit patterns. The absence of these forbidden digit patterns can be formulated with a CNF formula, 
 which then yields
an fpt-reduction to weighted CNF satisfiability with parameter $k$, a problem that is known to be in $\W[2]$.

We next turn to permutation groups, where addition is replaced by
composition of permutations. Given a list of permutations, a target
permutation, and a length bound $k\in\N$, the \emph{length problem} asks whether
the target is a product of at most $k$ permutations from the list.
Repetitions are allowed, and the order of the chosen permutations is unrestricted.
Even and Goldreich \cite{EvenGoldreich1981} proved that the problem is
$\mathsf{NP}$-complete when $k$ is encoded in unary.
For binary-encoded $k$, Jerrum \cite{Jerrum1985} proved
$\mathsf{PSPACE}$-completeness, even for two input generators.
In a recent preprint, Thumm and the author \cite[Theorem~1.2]{LohreyThumm2026}
established a sharp distinction for nilpotent groups: the binary length
problem is $\mathsf{NP}$-complete for $2$-step nilpotent permutation groups and
$\mathsf{PSPACE}$-complete already for permutation groups that are simultaneously
$3$-step nilpotent and metabelian.

When the length bound $k$ is the parameter, the exact-length version (where a product
of exactly $k$ permutations from the input list must yield the target)
 is called \textsc{Permutation Group Factorization} by Cai, Chen, Downey, and Fellows
\cite{CaiChenDowneyFellows1997}. They proved that it is $\W[1]$-hard
\cite[Theorem~6]{CaiChenDowneyFellows1997}. In the same paper, they proved
$\W[2]$-hardness for the corresponding factorization problem in
transformation monoids, whose elements are arbitrary self-maps of a finite
set \cite[Theorem~7]{CaiChenDowneyFellows1997}.
Lohrey and Rosowski \cite{LohreyRosowski2025} denote \textsc{Permutation Group
Factorization} by $\pF[S_*]$, where $\pF$ stands for parameterized
factorization and $S_*$ is the class of all symmetric groups.
In the present paper, we write $\Sym(*)$ for this class and hence
$\pF[\Sym(*)]$ for the problem. 
The related parameterized subset sum problem $\pSSS[\Sym(*)]$ asks whether a target permutation is an
ordered product of $k$ entries selected from a given list of permutations, with each
list element used at most once and the list order preserved.
There is a simple fpt-reduction from $\pF[\Sym(*)]$ to $\pSSS[\Sym(*)]$
that works for arbitrary groups and does not change the ambient group generated
by the input group elements; see \cref{pre:lem:replication}. 
In the rest of the introduction we only consider the exact-length variants $\pF[\Sym(*)]$
and $\pSSS[\Sym(*)]$, but all our results also hold for the variants where the target has to 
be produced by at most $k$ elements from the input list.

The best known upper bound for $\pSSS[\Sym(*)]$ is membership in
$\W^{\func}[2]$, established by Lohrey and Rosowski
\cite{LohreyRosowski2025}. The difference between $\W^{\func}[2]$ and $\W[2]$ can be expressed using model checking:
both classes can be defined as the fpt-closure of the model-checking problem for finite structures and
first-order sentences with an existential
quantifier block followed by one universal quantifier, but $\W[2]$
uses relational vocabularies, whereas $\W^{\func}[2]$ also permits
function symbols and nested terms. A machine-based characterization of $\W^{\func}[2]$ can be found in \cite{Chen2004}.
By definition, we have $\W[2]\subseteq\W^{\func}[2]$; the reverse inclusion is not known.
Despite the upper bound $\W^{\func}[2]$ for $\pSSS[\Sym(*)]$, it remains open whether
$\pSSS[\Sym(*)]$ belongs to any fixed level of the $\W$-hierarchy.

When restricting $\pSSS[\Sym(*)]$ (or $\pF[\Sym(*)]$) to the case where the input permutations
generate a solvable permutation group $G$ of derived length $d \ge 1$, where $d$ is a
fixed constant, we obtain an upper bound of $\W[2d]$
(\cref{thm:solvable}). 
The proof is by induction on $d$ and constructs from the input instance an instance of 
the weighted satisfiability problem for Boolean circuits of weft at most $2d$ and depth at most $4d$.
We proceed along the derived series of $G$ and equip the orbits of
the abelian normal subgroup $G^{(d-1)} \unlhd G$ with abelian group structures.
An affine decomposition expresses product evaluation in $G$ through an action of a group 
of smaller derived length and sums in these orbit groups (\cref{lem:affine-product}).
Encoding the sums by positive CNFs (see \cref{lem:abelian-cnf}) adds at most two to the circuit weft
and four to its depth at each step, giving circuits of weft at most
$2d$ and depth at most $4d$.

For nilpotent permutation groups of any fixed nilpotency class, we can further lower the complexity 
to $\W[2]$ (\cref{thm:nilpotent}). The main ingredients are Hall coordinates and the Magnus embedding to express an ordered
product equation in a nilpotent group by polynomial congruences of bounded degree.
An fpt-reduction then replaces these congruences by a subset sum instance
over an abelian permutation group, which allows us to use the above-mentioned \cref{thm:solvable} for $d=1$.

Finally, we improve the lower bound for $\pF[\Sym(*)]$ (and hence also for $\pSSS[\Sym(*)]$) from $\W[1]$ to 
$\W[2]$ (\cref{thm:nilpotent-hardness}). This lower bound holds already for the restriction where the input permutations generate a $2$-group (a group in which every element has order $2^p$ for some $p$), which is necessarily nilpotent.
Our reduction starts from parameterized \textsc{Set Cover} (given sets $S_1, \ldots, S_m \subseteq U$, are there at most $k$ subsets $S_i$ that cover $U$?)
and uses nested commutators in a dihedral $2$-group to test, for each
element $a \in U$ of the \textsc{Set Cover} instance, whether at least one selected set $S_i \subseteq U$ contains $a$.
Additional cyclic factors of order two ensure that these selected $S_i$, when taking all $a \in U$ into account, form
a consistent choice of at most $k$ sets from $S_1, \ldots, S_m$.
Finally, copies of the dihedral group of order eight force the group elements
to occur in the required order that is imposed by the constructed $\pF[\Sym(*)]$ instance.
The nilpotency class of the constructed groups depends on the parameter $k$ and is therefore not fixed.

\section{Preliminaries}

In this section we collect general notations and define the necessary framework 
from parameterized complexity and group theory (in particular, permutation groups).

\subsection{General notations}

We write $\N=\{0,1,2,\ldots\}$. For integers $i,j$
we write $[i,j]$ for the interval $\{i,i+1,\ldots,j\}$; if $i>j$, this
interval is empty. For a Boolean condition $A$ that is either true or false
we define $\mathbf{1}[A]$ as $1$ if $A$ is true and otherwise as $0$.
We write $\poly(n)$ for an arbitrary polynomial in $n$.

For a row vector $\overline{a} = (a_1, a_2, \ldots, a_k) \in \mathbb{R}^k$ we define the norm
$|\overline{a}|_1 = |a_1|+|a_2|+\cdots +|a_k|$. 
The least common multiple of integers $a_1, \dots, a_k$
 is denoted by $\operatorname{lcm}\{ a_1, \ldots, a_k\}$.
With $\overline{0}$ we denote a tuple of zeros, whose length is determined by the context and similarly for $\overline{1}$.
For a row vector $\overline{a}$ we write $\overline{a}^{\mathsf{T}}$ for the transposed column vector.

For a multivariate polynomial $P \in \mathbb{R}[\overline{x}]$ we define the (total) \emph{degree} of $P$ 
as follows: For a monomial $x_1^{n_1} x_2^{n_2} \cdots x_t^{n_t}$ the degree is $n_1 + n_2 + \cdots + n_t$
and the degree of $P$ is the maximal degree of a monomial in $P$. For the zero polynomial we define the degree as $-\infty$.

For a finite alphabet $\Sigma$ of symbols, we denote by $\Sigma^*$ the set of all finite words
over $\Sigma$. The length of a word $w \in \Sigma^*$ is denoted by $|w|$. 
For $n \in \N$ we define $\Sigma^n = \{ w \in \Sigma^*:|w|=n\}$ and $\Sigma^{\le n} = \{ w \in \Sigma^*:|w| \le n\}$.

\subsection{Parameterized complexity}

For a detailed introduction to parameterized complexity theory, the reader should consult the monograph \cite{FlumG06}.
A \emph{parameterized problem} is a subset
$L\subseteq\Sigma^*\times\N$ for a finite alphabet $\Sigma$. For a pair
$(u,k)\in\Sigma^*\times\N$, we call $k$ the parameter. For a fixed number
$d\in\N$, the $d$-th slice of $L$ is
$L_d=\{(u,d)\in L:u\in\Sigma^*\}$. The class $\FPT$ of
\emph{fixed-parameter tractable problems} consists of all parameterized
problems that can be decided in time $f(k) \cdot \poly(|u|)$ for a computable
function $f:\N\to\N$.

An \emph{fpt-reduction} from $L_1\subseteq\Sigma_1^*\times\N$ to
$L_2\subseteq\Sigma_2^*\times\N$ is a computable function $R : \Sigma_1^*\times\N \to \Sigma_2^*\times\N$ 
 such that the following hold for all $(u,k) \in \Sigma_1^*\times\N$:
\begin{itemize}
\item $(u,k)\in L_1$ if and only if $R(u,k)\in L_2$;
\item $R(u,k)$ can be computed in time $f(k)\cdot\poly(|u|)$ for a computable
function $f:\N\to\N$;
\item there is a computable function $h:\N\to\N$ such that $R(u,k)=(u',k')$ implies $k'\leq h(k)$.
\end{itemize}
We write $L_1\lfpt L_2$ if there is an fpt-reduction from $L_1$ to $L_2$, and
$L_1\equiv^{\fpt}L_2$ if $L_1 \lfpt L_2 \lfpt L_1$. For a
parameterized problem $L$, its fpt-closure $[L]^{\fpt}$ is the class of all
parameterized problems $K$ with $K\lfpt L$.

We use the circuit characterization of the $\W$-hierarchy
\cite[Remark~3.3 and Section~4.3]{FlumGrohe2005}. Boolean circuits are built up from
and- and or-gates with arbitrary fan-in and input gates that are labelled with a Boolean variable or 
a negated Boolean variable. An and-gate (resp., or-gate) with fan-in $0$ can be identified with a constant gate
carrying the truth value $1$ (resp., $0$).
A gate is \emph{large} if its fan-in is at least three. The \emph{weft} of a
Boolean circuit is the maximum number of large gates on a path
in the circuit, and its depth is the maximum number of
wires on a path. If $X$ is the set of Boolean variables that appear (possibly negated) at
the input gates of a circuit then a \emph{truth assignment} for the circuit is a mapping
$\phi : X \to \{0,1\}$. Under such a truth assignment, every circuit gate $B$ evaluates to a truth
value $\phi(B)$ in the obvious way. The \emph{weight} of a truth assignment $\phi  : X \to \{0,1\}$ is 
the number of $x \in X$ with $\phi(x)=1$.
For fixed $d,c\ge 1$, we define the following problem
\begin{problem}[{weighted circuit satisfiability for depth $d$ and weft $c$, $\WCS[d,c]$ for short}]\rm~\\
Input: a Boolean circuit $C$ with a distinguished output gate $B$, depth at most $d$ and weft at most $c$ and $k\ge 1$.\\
Parameter: $k$.\\
Question: Is there a truth assignment $\phi$ for $C$ of weight $k$ such that $\phi(B)=1$?
\label{problem:WCS}
\end{problem}
\noindent
The class $\W[k]$ for $k \geq 1$ is defined as $\W[k] = \bigcup_{d \geq 1}[\WCS[d,k]]^{\fpt}$.

Recall that a Boolean formula is in \emph{conjunctive normal form} (CNF for short) if it is a conjunction of clauses,
where a clause is a disjunction of literals (Boolean variables or negated Boolean variables).
A CNF is \emph{positive} if all its literals are variables. As usual, we identify a clause with a set of literals.
The problem $\wCNFSAT$ is defined analogously to $\WCS[d,c]$ but the Boolean circuit $C$ is replaced by a CNF. 
This problem is $\W[2]$-complete \cite[Theorem~4.13]{FlumGrohe2005}.

\subsection{Groups}

We assume that the reader is familiar with basic concepts in group theory (like normal subgroups, quotient groups, group homomorphisms, etc.);
see \cite{Robinson1996} for more details. We write $H \unlhd G$ if $H$ is a normal subgroup of $G$. 
For a group $G$, an element $g \in G$ is \emph{central} if $g h = hg$ for all $h \in G$.
For elements $g, h \in G$ we define the \emph{commutator} $[g,h]=g^{-1}h^{-1}gh$. 
For $A, B \subseteq G$ we write $[A,B]$ for the subgroup of $G$ generated by all commutators
$[a,b]$ with $a \in A$ and $b \in B$.
We define the \emph{lower central series} of $G$ by $\gamma_1(G)=G$ and
$\gamma_{i+1}(G)=[\gamma_i(G),G]$ for $i\ge 1$. Every $\gamma_i(G)$ is a normal subgroup of $G$ and $\gamma_{i+1}(G) \unlhd \gamma_i(G)$.
Moreover, $[\gamma_i(G),\gamma_j(G)]\subseteq\gamma_{i+j}(G)$ holds for all $i,j \geq 1$.
A group $G$ is called \emph{$d$-step nilpotent} if $\gamma_{d+1}(G)=1$ and $G$ is called \emph{nilpotent} if $G$ is $d$-step nilpotent for some $d$.
For a prime $p$, a finite $p$-group is a finite group whose order is a power of $p$.
Every finite $p$-group is nilpotent \cite[5.1.3]{Robinson1996}.

The \emph{derived series} of a group $G$ is defined by $G^{(0)}=G$ and
$G^{(i+1)}=[G^{(i)},G^{(i)}]$ for $i\in\N$. Every $G^{(i)}$ is a normal subgroup of $G$.
If $G^{(d)}=1$ then $G$ is called \emph{solvable of derived length $d$} and if there is a $d$ with 
$G^{(d)}=1$ then $G$ is called \emph{solvable}.
Note that our definition implies that if $G$ is solvable of derived length $d$ then it is also solvable of derived length $d+1$.
 
 For a group $G$ and a subset $A \subseteq G$ we write $\langle A \rangle$ for the smallest subgroup of
 $G$ containing $A$. If $G = \langle A \rangle$ then $A$ is called a \emph{generating set} for $G$.
 If $A = \{g_1, \ldots, g_m\}$ then we also write $\langle g_1, \ldots, g_m\rangle$ for $\langle A \rangle$.
 
\subsection{Permutation groups} \label{sec:PGs}

For a set $\Omega$ we denote with $\Sym(\Omega)$ the group of all permutations on the set $\Omega$.
In this paper, $\Omega$ will always be finite. If $\Omega=[1,n]$ then we also write $\Sym(n)$ for $\Sym(\Omega)$. 
For $g \in \Sym(n)$  and $a \in [1,n]$ we write $a^g$ instead of $g(a)$. As a consequence, we have
$x^{gh}=(x^g)^h$ for $x\in [1,n]$ and $g,h\in \Sym(n)$.
A (finite) \emph{permutation group} $G$ is a subgroup of some $\Sym(n)$; $n$ is then also called the \emph{permutation degree} of $G$.\footnote{The permutation degree
of a permutation group is not an intrinsic property of the abstract group $G$ but depends on the concrete embedding of $G$ into a symmetric group.}
The description of a permutation group $G \le \Sym(n)$ is a list $g_1, \ldots, g_m \in \Sym(n)$
with $G = \langle g_1, \ldots, g_m\rangle$. 
Each permutation $g \in \Sym(n)$ is represented by
 its \emph{one-line representation}, which is the list $[1^g, 2^g, \ldots, n^g]$.

By Cayley's theorem every finite group $G$ of size $n$ is isomorphic to a subgroup of $\Sym(n)$.
The important point about permutation groups is that the description of a permutation group by a set
of generating permutations can be much more succinct than the description by a multiplication table.
For instance, $\Sym(n)$ has size $n!$ but it can be generated by two permutations.

In this paper, we will use a few standard algorithms for permutation groups; see 
 \cite{Seress2003} for more details.  
 Given a permutation group $G \leq \Sym(n)$ and an element
 $g \in \Sym(n)$ one can check in polynomial time whether $g \in G$ and compute $|G|$.
 By a straightforward enumeration, one can compute a list of all $G$-orbits $\{a^g : g \in G\}$ ($a \in [1,n]$)
 in polynomial time. 
 Finally, the derived subgroup $G^{(1)}$ of a given permutation group $G$
can also be computed in polynomial
time \cite{FurstHopcroftLuks1980}: Start with the permutation group generated by all commutators $[a,b]$, where $a$ and $b$ are
generators of $G$ and close this group under conjugation by generators of $G$.
Add a conjugate only if membership testing shows
that it enlarges the current subgroup. There are at most
$\log_2|G|\leq\log_2(n!)$ enlargements of the group. The resulting normal closure
is $G^{(1)}$, since the images of the original generators commute in the
quotient. This construction can be iterated to compute the derived series of a given permutation group $G$
in polynomial time. In the same way, one can also compute the lower central series of 
a given permutation group $G$ in polynomial time. In particular, one can check in polynomial time whether a given
permutation group is nilpotent, $d$-step nilpotent (for a fixed $d$), solvable, or solvable of derived length $d$ (for a fixed $d$).

\subsection{Straight-line programs}

A \emph{straight-line program} (SLP) over an alphabet $\Sigma \neq \emptyset$ 
is a sequence of instructions 
\begin{equation} \label{eq:SLP}
x_1 := \alpha_1; x_2 := \alpha_2; \cdots x_m := \alpha_m,
\end{equation}
where every $\alpha_i$ is either a symbol from $\Sigma$ or the concatenation $x_j x_k$ of two 
variables $x_j,x_k$ with $j,k < i$. The size of this SLP is $m$. An SLP can be evaluated
in a monoid $M$ provided every alphabet symbol $a \in \Sigma$ is interpreted as an element of $M$.
Given a finite nonempty subset $A \subseteq M$ and an SLP $\mathcal{P}$ over the alphabet $A$, we say that
$\mathcal{P}$ is an SLP for $a \in M$ if the value of the last variable of $\mathcal{P}$ ($x_m$ in \ref{eq:SLP})
is $a$.

By evaluating an SLP over an alphabet $\Sigma$ in the free monoid $\Sigma^*$, one obtains
a word, whose length can be exponential in the size of the SLP \cite{Lohrey12}. SLPs are frequently used in the
context of black-box groups \cite{BabaiS84}. For permutation groups, we use the following result from \cite{FurstHopcroftLuks1980}.

\begin{theorem} \label{theorem:slp}
Given $n \geq 1$, a finite nonempty set $A \subseteq \Sym(n)$ and $g \in \Sym(n)$, one can check in polynomial time whether 
$g \in \langle A \rangle$, and in case this holds, one can compute in polynomial time an SLP for $g$ over $A$
of size $\poly(n)$.
\end{theorem}

\section{Subset sum and factorization for classes of groups}
\label{sec:generic-problems}

In order to obtain a uniform description of the problems considered in this paper, we use the framework from \cite[Section~4.1]{LohreyRosowski2025}.
Let $\mathsf C$ be a countable class of countable groups.
We assume that there exists a description scheme that associates to every group in $\mathsf C$
a unique finite description (e.g., a bit string) that will be part of the input below. In addition, also 
elements of groups in $\mathsf C$ have a finite description.
The description scheme of $\mathsf C$ is part of the problem specification. In
particular, $\mathsf C$ may have two different description schemes that lead to different
computational problems. 
In all group descriptions used below, group multiplication, comparison, and construction of the
group identity are polynomial-time operations.

For a group $G$, a finite subset $A\subseteq G$ and $k\in\N$, let
\[
 A^k=\{a_1a_2\cdots a_k:a_1,a_2,\ldots,a_k\in A\} \subseteq G,
 \qquad A^{\leq k}=\bigcup_{s=0}^k A^s,
\]
where $A^0=\{1\}$.  We will consider the following parameterized problems for a class of groups $\mathsf C$ with the above properties.

\begin{problem}[{parameterized subset sum for $\mathsf C$, $\pSSS[\mathsf C]$ for short}]\rm~\\
Input: the description of a group $G\in\mathsf C$, a list
$h, g_1,\ldots,g_m$ of elements from $G$, and $k\ge 1$.\\
Parameter: $k$.\\
Question: Are there $\epsilon_1,\ldots,\epsilon_m\in\{0,1\}$ such that
$\sum_{i=1}^m\epsilon_i=k$ and $h = g_1^{\epsilon_1}\cdots g_m^{\epsilon_m}$ hold?
\label{problem:SSS}
\end{problem}
Note that $g_i = g_j$ for $i < j$ may hold.
An instance with $k>m$ is necessarily negative. We therefore assume throughout this paper that $k \in [1,m]$ in instances
of $\pSSS[\mathsf C]$.

\begin{problem}[{parameterized factorization for $\mathsf C$, $\pF[\mathsf C]$ for short}]\rm~\\
Input: the description of a group $G\in\mathsf C$, elements
$h,g_1,\ldots,g_m \in G$, and $k\ge 1$.\\
Parameter: $k$.\\
Question: Is $h\in\{g_1,\ldots,g_m\}^k$?
\label{problem:F}
\end{problem}
In contrast to $\pSSS[\mathsf C]$, in $\pF[\mathsf C]$, both the order and the multiplicities of the
group elements $g_i$ may be chosen. 
We write $\pSSS_{\leq}[\mathsf C]$ and $\pF_{\leq}[\mathsf C]$ for the
variants with $\sum_{i=1}^m\epsilon_i \le k$ in $\pSSS[\mathsf C]$
and $h\in\{g_1,\ldots,g_m\}^{\le k}$ in $\pF[\mathsf C]$.
For abelian groups we use additive notation; hence, we write $\sum_{i=1}^m \epsilon_i \cdot g_i=h$ in 
$\pSSS_{\leq}[\mathsf C]$.
The following simple fpt-reductions can be found in \cite[Lemmas~4.7 and 4.16]{LohreyRosowski2025}.

\begin{lemma}
\label{pre:lem:replication}
For every class $\mathsf C$ with the representations just described, we have
$\pF[\mathsf C]\lfpt\pSSS[\mathsf C]$ and 
$\pF_{\leq}[\mathsf C]\lfpt\pSSS_{\leq}[\mathsf C]$.
Both reductions preserve the group $\langle g_1, \ldots, g_m\rangle$ and the parameter $k$.
\end{lemma}
For the proof of \cref{pre:lem:replication}, one only has to repeat the list $g_1, \ldots, g_m$ of input group elements $k$ times.
By adding $k$ copies of the group identity $1$ to the input list $g_1, \ldots, g_m$, we obtain the following lemma:
\begin{lemma}
\label{pre:lem:padding}
For every class $\mathsf C$, for which the group identity of a given $G \in \mathsf C$ can be constructed
in polynomial time, we have
$\pSSS_{\leq}[\mathsf C]\lfpt\pSSS[\mathsf C]$ and $\pF_{\leq}[\mathsf C]\lfpt\pF[\mathsf C]$.
Both reductions preserve the group $\langle g_1, \ldots, g_m\rangle$ and the parameter $k$.
\end{lemma}
In \cref{sec:integers} we consider the classical subset sum problem over the integers
with binary-encoded integers. It is obtained by taking $\mathsf C = \{ \Z\}$ with the 
convention that elements of $\Z$ are given in binary notation.
In Sections~\ref{sec:solvable}--\ref{sec:nilpotent-hardness} we deal with 
the class $\mathsf C = \{ \Sym(n) : n \geq 1 \}$ of symmetric groups on a finite domain. 
The representation of a symmetric group $\Sym(n)$ simply consists of the number $n$ (it can be encoded in unary notation).
Let us write $\Sym(*)$ for this class. Later we will consider restrictions of the problem $\pSSS[\Sym(*)]$ and its variants, where 
 the permutation group generated by the input permutations $\langle g_1, \ldots, g_m \rangle$ has to satisfy certain algebraic properties.

\section{Subset sum for integers}
\label{sec:integers}

In this section, we will show that $\pSSS[\Z]$ belongs to $\W[2]$, which improves the upper bound of $\W[3]$ from \cite{BussIslam2007}.
The crucial technical contribution is \cref{lem:integer-carries} that allows us to reduce $\pSSS[\Z]$ to $\wCNFSAT$.

For $b,\ell\in\N$ with $b\geq2$ and $\ell\geq1$, a number $x \in [0,b^\ell-1]$, and $p \in [0,\ell-1]$, we write $\dig_b(x,p)$ for the digit of $x$ at position $p$
when $x$ is written in base $b$. Formally, the digits $\dig_b(x,p)\in[0,b-1]$ are uniquely defined by the condition
$x = \sum_{p \in [0,\ell-1]} \dig_b(x,p) \cdot b^p$. For $p\in\N$ with $p \geq \ell$, we define $\dig_b(x,p)=0$.

In the following lemma we consider the addition of $k$ numbers $a_1, \ldots, a_k \in [0,(k+1)^\ell-1]$ each having
at most $\ell$ digits when written in base $k+1$. The digit positions of these numbers are numbered from $0$ to $\ell-1$.
\begin{lemma}
\label{lem:integer-carries}
Let $k,\ell\ge 1$ and
$a_1,\ldots,a_k,t\in[0,(k+1)^\ell-1]$.
For $i \in [1,k]$ and $p \in [0,\ell]$ define
\begin{itemize}
\item $a_{i,p} = \dig_{k+1}(a_i,p)$, 
\item $t_p =  \dig_{k+1}(t,p)$,
\item $s_p=\sum_{i=1}^k a_{i,p}\in\N$ and
\item $c_p=(t_p-s_p)\bmod (k+1)\in[0,k]$.
\end{itemize}
Then the following conditions are equivalent:
\begin{enumerate}[\rm (1)]
\item $\sum_{j=1}^k a_j=t$.
\label{eq:integer-local-tests-1}
\item $c_0=0$ and $\forall p\in[0,\ell-1]: s_p+c_p=t_p+(k+1) \cdot c_{p+1}$.
 \label{eq:integer-local-tests}
\end{enumerate}
\end{lemma}

\begin{proof}
Note that $a_{i,\ell} = t_\ell = s_\ell = c_\ell = 0$.
Suppose first that \eqref{eq:integer-local-tests-1} holds.
Let $c'_p\in\N$ be the carry that arrives in position $p\in[0,\ell]$ when adding the $a_i$ in base $k+1$.
We claim that $c'_p \in [0,k-1]$. Indeed, $c'_\ell = c'_0 = 0$ and if $c'_p \in [0,k-1]$ then 
$c'_{p+1} \leq \lfloor(k^2+(k-1))/(k+1)\rfloor=k-1$. 

Since $c'_p \in [0,k-1]$, we have $t_p = (s_p + c'_p) \bmod (k+1)$, i.e.,  $c'_p = (t_p-s_p) \bmod (k+1)$.\footnote{For this argument, $c'_p \in [0,k]$ would suffice.}
We therefore have $c'_p = c_p$. Since  $s_p+c'_p=t_p+(k+1) \cdot c'_{p+1}$ we obtain
\eqref{eq:integer-local-tests}.

Conversely, \eqref{eq:integer-local-tests} implies
\begin{align*}
 \sum_{j=1}^k a_j-t
 =\sum_{p=0}^{\ell-1}(s_p-t_p)(k+1)^p=\sum_{p=0}^{\ell-1}((k+1)c_{p+1}-c_p)(k+1)^p
 =(k+1)^\ell c_\ell-c_0=0,
\end{align*}
where the last equality follows from $c_0=0$ and $c_\ell = (t_\ell-s_\ell)\bmod (k+1) = 0$.
\end{proof}

\begin{example}
In \cref{lem:integer-carries}, let $k=\ell=3$ and choose
$a_1,a_2,a_3,t\in[0,4^3-1]$ with
$a_1=27$, $a_2=10$, $a_3=7$, and $t=44$ in decimal notation.
Their base-four representations are
$a_1=(123)_4$, $a_2=(022)_4$, $a_3=(013)_4$, and $t=(230)_4$.
Using the notation of the lemma, we obtain the following table:
\begin{center}
\begin{tabular}{c|rrrr}
$p$       & $0$ & $1$ & $2$ & $3$ \\ \hline
$a_{1,p}$ & $3$ & $2$ & $1$ & $0$ \\
$a_{2,p}$ & $2$ & $2$ & $0$ & $0$ \\
$a_{3,p}$ & $3$ & $1$ & $0$ & $0$ \\ \hline
$s_p$     & $8$ & $5$ & $1$ & $0$ \\
$t_p$     & $0$ & $3$ & $2$ & $0$ \\
$c_p$     & $0$ & $2$ & $1$ & $0$
\end{tabular}
\end{center}
For instance, $c_1=(t_1-s_1)\bmod 4=(3-5)\bmod 4=2$.
The equations $s_p+c_p=t_p+4c_{p+1}$ for $p=0,1,2$
are respectively $8+0=0+4\cdot2$, $5+2=3+4\cdot1$,
and $1+1=2+4\cdot0$.
Moreover, $c_0=0$. Thus condition~\eqref{eq:integer-local-tests}
holds, as expected from $27+10+7=44$.

The condition $c_0=0$ is essential. If we replace the target by
$t=45=(231)_4$, then $c_0=1$, $c_1=2$, $c_2=1$, and $c_3=0$.
All three local equations still hold: the first becomes
$8+1=1+4\cdot2$, and the other two remain unchanged.
However, $c_0 \neq 0$, so condition~\eqref{eq:integer-local-tests}
fails.
\end{example}

\begin{theorem}
\label{thm:integers}
The problem $\pSSS[\Z]$ belongs to $\W[2]$.
\end{theorem}

\begin{proof}
We show that $\pSSS[\Z] \lfpt \wCNFSAT$.
Consider an instance $a_1,\ldots,a_m,t\in\Z$, $k\in [1,m]$ of $\pSSS[\Z]$. 
Choose $u\in\N$ with
$a_i+u\geq0$ for all $i\in[1,m]$. Replacing $a_i$ by $a_i+u$ and
$t$ by $t+ku$ preserves all solutions of cardinality $k$. If the new
target $t$ is negative, the instance is negative. We may therefore assume
that $a_1,\ldots,a_m,t\in\N$.

Choose $\ell\ge 1$
such that all $a_1, \ldots, a_m, t  < (k+1)^\ell$. For $i \in [1,m]$ and $p \in [0,\ell-1]$ let
 $a_{i,p} = \dig_{k+1}(a_i,p)$, $t_p =  \dig_{k+1}(t,p)$, and $a_{i,\ell} = t_\ell = 0$.

We now define an instance of $\wCNFSAT$. The parameter will be $k$.
Introduce propositional variables $x_{r,i}$ for $r\in[1,k]$ and
$i\in[1,m]$. The intended meaning is that the $r$-th selected index from $[1,m]$ is $i$. For every $r\in[1,k]$, include the clause
\begin{equation} \label{eq:selection-clause}
C_r = \{ x_{r,i} : 1 \le i \le m \}.
\end{equation} 
Since we ask for a satisfying assignment of weight exactly
$k$, the clauses $C_r$ enforce that, for every $r\in[1,k]$, exactly one variable $x_{r,i_r}$ is true.
For all $r,s\in[1,k]$ with $r<s$ and $i\in[1,m]$, also include the clause
\begin{equation} \label{eq:selection-clause2}
C_{r,s,i} = \{ \neg x_{r,i}, \neg x_{s,i} \}.
\end{equation} 
These clauses prevent an index $i \in [1,m]$ from
being chosen twice.

For $p\in[0,\ell-1]$, consider tuples
$\bar{\alpha} = (\alpha_1,\ldots,\alpha_k), \bar{\alpha}' = (\alpha'_1,\ldots,\alpha'_k)\in [0,k]^k$. They describe
proposed digits of the selected numbers at positions $p$ and $p+1$.
Define $a=\sum_{r=1}^k\alpha_r$,
$a'=\sum_{r=1}^k\alpha'_r$,
$c=(t_p-a)\bmod (k+1)\in [0,k]$, and
$c'=(t_{p+1}-a')\bmod (k+1)\in [0,k]$.
Call the pair $(\bar{\alpha},\bar{\alpha}')$ \emph{admissible at $p$} if
$a+c=t_p+(k+1)c'$ and, when $p=0$, also $c=0$.
For every pair $(\bar{\alpha},\bar{\alpha}')$ with $\bar{\alpha} = (\alpha_1,\ldots,\alpha_k), \bar{\alpha}' = (\alpha'_1,\ldots,\alpha'_k)$
 that is not admissible at $p$, include the clause
\begin{equation}
C_{\bar{\alpha},\bar{\alpha}',p}  =  \{ x_{r,i} 
:  1 \leq r \leq k, 1 \leq i \leq m, (a_{i,p},a_{i,p+1})\neq(\alpha_r,\alpha'_r) \}.
 \label{eq:forbidden-clause}
\end{equation}
As usual, an empty clause is false.

Let $i_1,\ldots,i_k\in[1,m]$ be the distinct indices determined by an assignment $\phi$ of weight $k$ satisfying all clauses from \eqref{eq:selection-clause}
and \eqref{eq:selection-clause2}. Hence, for every $r \in [1,k]$, $\phi(x_{r,i_r})=1$ and
$\phi(x_{r,i})=0$ for every $i\in[1,m]\setminus\{i_r\}$.
Then the clause $C_{\bar{\alpha},\bar{\alpha}',p}$ is false under the assignment $\phi$ if and only if
$(a_{i_r,p},a_{i_r,p+1})=(\alpha_r,\alpha'_r)$ for all $r \in [1,k]$, which means that
$((a_{i_1,p}, \ldots, a_{i_k,p}), (a_{i_1,p+1}, \ldots, a_{i_k,p+1})) = (\bar{\alpha},\bar{\alpha}')$, where by construction
$(\bar{\alpha},\bar{\alpha}')$ is not admissible at $p$.
Thus, the assignment $\phi$ satisfies all clauses in \eqref{eq:forbidden-clause} if and only if for every $p \in [0,\ell-1]$ the pair
$((a_{i_1,p}, \ldots, a_{i_k,p}), (a_{i_1,p+1}, \ldots, a_{i_k,p+1}))$ is admissible at $p$.
 By \cref{lem:integer-carries}, this is equivalent to
$\sum_{r=1}^k a_{i_r}=t$. This proves correctness in both directions.

The above construction is an fpt-reduction: there are $km$ variables $x_{r,i}$ and the
number of clauses is bounded by $k+\binom{k}{2}m+\ell (k+1)^{2k}$, where $\ell$ is linearly bounded in the maximal bit length
of the input numbers $a_1, \ldots, a_m, t$.
Moreover, the parameter $k$ does not change.
\end{proof}
With \cref{pre:lem:padding,pre:lem:replication} we obtain:
\begin{corollary}
\label{coro:integers}
The problems $\pF_{\le}[\Z]$, $\pF[\Z]$, $\pSSS_{\le}[\Z]$ and $\pSSS[\Z]$ belong to $\W[2]$.
\end{corollary}

\section{Subset sum for solvable permutation groups}
\label{sec:solvable}

In this section, $d \ge 1$ will be a fixed constant. 
We then consider the restriction of the problem $\pSSS[\Sym(*)]$, where the input permutations $g_1, \ldots, g_m$ generate a
solvable permutation group of derived length $d$. Recall from \cref{sec:PGs} that this property can be checked in polynomial time.
We denote this restriction of $\pSSS[\Sym(*)]$ with $\pSSS[\SolPG{d}]$; the problems $\pSSS_\le[\SolPG{d}]$, $\pF[\SolPG{d}]$ and
$\pF_\le[\SolPG{d}]$ are defined similarly.

\subsection{Positive CNFs for addition in small abelian groups}
\label{sec:abelian-cnf}

In this section we consider finite abelian groups that are given by their multiplication table.
The following standard algorithmic fact also follows from the deterministic
basis algorithm of Chen and Fu \cite{ChenFu2011}: once a basis is known,
enumerating its coordinate tuples gives the full table of the isomorphism.
We include a proof using Smith normal form.

\begin{lemma}
\label{lem:abelian-compute-cyclic-deco}
Let $H$ be a finite abelian group
given by its addition table. In polynomial time one can compute integers
$q_1, \ldots, q_m \in [2,|H|]$ and an isomorphism $\theta : H \to \prod_{1 \leq i \leq m} \Z/q_i\Z$, represented by a table containing all values $\theta(a)$ for $a \in H$.
\end{lemma}

\begin{proof}
We form the integer matrix $A$ with columns $e_x+e_y-e_{x+y}$ for $x,y\in H$,
where $e_x\in\Z^H$ is the corresponding standard basis vector.
These relations present $H$. The Smith normal form $S$ of $A$ yields the decomposition
of $H$ into cyclic subgroups \cite{Iliopoulos1989}. The matrix $S$ together with unimodular integer matrices $U$ and $V$
such that $A = USV$ can be computed in polynomial time \cite{KannanBachem1979}. 
The moduli $q_i$ are the diagonal entries of $S$ that are different from $0$, $1$ and $-1$.
To transform $e_x$ into the corresponding tuple $\theta(x)$ one computes
$U^{-1} e_x$ and reduces the entries modulo the $q_i$.
\end{proof}

\begin{lemma}
\label{lem:abelian-cnf}
Let $k\ge 1$, let $H$ be a finite abelian group
given by its addition table, and put $n=|H|\ge 1$.
For $h\in H$, one can compute a positive CNF $\Phi[H,k,h]$
over the propositional variables $y_{r,a}$, where $r\in[1,k]$ and $a\in H$,
with the following property. For all $h_1,\ldots,h_k\in H$ and the corresponding
truth assignment $y_{r,a} := \mathbf{1}[a=h_r]$, we have
\begin{equation}
 \Phi[H,k,h]=1
 \quad\Longleftrightarrow\quad
 \sum_{r=1}^{k}h_r=h\text{ in }H.
 \label{eq:small-abelian-cnf}
\end{equation}
The size and construction time are bounded by $f(k) \cdot n^c$ for a computable
function $f:\N\to\N$ and a constant $c\in\N$ independent of $k$.
\end{lemma}

\begin{proof}
Using \cref{lem:abelian-compute-cyclic-deco},
we first compute in polynomial time a cyclic decomposition $\prod_{1 \leq i \leq m} \Z/q_i\Z$ 
and an isomorphism $\theta : H \to \prod_{1 \leq i \leq m} \Z/q_i\Z$, represented by a table containing all values $\theta(a)$ for $a \in H$.
For $a \in H$ we write 
\[ \theta(a) = (\theta_1(a), \ldots, \theta_m(a)) \in \prod_{1 \leq i \leq m} \Z/q_i\Z = \prod_{1 \leq i \leq m} [0,q_i-1].\]
Fix an $i \in [1,m]$ for the following discussion. For $h_1,\ldots,h_k\in H$ we have
\begin{equation}
 \sum_{r=1}^{k}\theta_i(h_r)\equiv \theta_i(h) \pmod {q_i}
 \quad\Longleftrightarrow\quad
 \bigvee_{j=0}^{k-1}
 \left(\sum_{r=1}^{k}\theta_i(h_r)=\theta_i(h)+jq_i\right).
 \label{eq:small-modulus-quotients}
\end{equation}
Indeed, the sum of the $\theta_i(h_r)$ in $\Z$ is strictly smaller than $k q_i$.

Let $\ell\ge 1$ be minimal with $(k+1)^\ell\geq kq_i$. All $\theta_i(h_r)$ and all targets
$\theta_i(h)+jq_i$, for $j\in[0,k-1]$, are smaller than $kq_i$, and
$\ell\leq2+\log_2 n$. For each $j\in[0,k-1]$, apply the
construction from the proof of
\cref{thm:integers} to the equality $\sum_{r=1}^{k}\theta_i(h_r)=\theta_i(h)+jq_i$
in \eqref{eq:small-modulus-quotients}. The variables $x_{r,i}$
are replaced by the variables $y_{r,a}$ ($1 \le r \leq k$, $a \in H$). 
We use only the clauses corresponding to \eqref{eq:forbidden-clause}. The clauses in 
\eqref{eq:selection-clause} are not needed since we only consider
assignments that make for every $r \in [1,k]$ exactly one of the variables $\{ y_{r,a} : a \in H \}$ true
and the clauses in \eqref{eq:selection-clause2} are not needed since we allow $h_r = h_s$ for $r < s$
in the lemma.

More explicitly, for each digit position $p\in[0,\ell-1]$, enumerate
all pairs $(\bar{\alpha},\bar{\alpha}')$ with $\bar{\alpha} = (\alpha_1,\ldots,\alpha_k), \bar{\alpha}' = (\alpha'_1,\ldots,\alpha'_k)$
that are not admissible at $p$ (at most $(k+1)^{2k}$ many). Here, in the definition of ``admissible'' 
in the proof of \cref{thm:integers}, we have to replace the target digits $t_p = \dig_{k+1}(t,p)$ by the target digits
$\dig_{k+1}(\theta_i(h)+jq_i,p)$. Then we construct all clauses
\begin{equation*}
C_{\bar{\alpha},\bar{\alpha}',p}  =  \{ y_{r,a} 
:  1 \leq r \leq k, a \in H, (\dig_{k+1}(\theta_i(a),p), \dig_{k+1}(\theta_i(a),p+1)) \neq(\alpha_r,\alpha'_r) \}.
 \label{eq:forbidden-clause-H}
\end{equation*}
Let $\mathcal D_{i,j}$ be the set of all these clauses.
The conjunction 
\[ 
\Psi_{i,j} = \bigwedge_{C \in \mathcal D_{i,j}} C 
\] 
(for $i \in [1,m]$ and
$j\in[0,k-1]$) is a positive CNF 
with at most $\ell (k+1)^{2k}$ clauses and expresses the integer equality
$\sum_{r=1}^{k}\theta_i(h_r)=\theta_i(h)+jq_i$.
Distributivity gives
\begin{equation}
 \Psi_{i} := \bigvee_{j=0}^{k-1}\Psi_{i,j}
 \quad\equiv \quad
 \bigwedge_{C_0\in \mathcal D_{i,0}} \cdots \bigwedge_{C_{k-1}\in \mathcal D_{i,k-1}} \!\!\!
        (C_0\vee\cdots\vee C_{k-1}).
 \label{eq:quotient-cnf-distribution}
\end{equation}
Each clause on the right of \eqref{eq:quotient-cnf-distribution}
is the union of $k$ positive clauses. There are at most
$(\ell (k+1)^{2k})^k$ such clauses. This is an fpt-bound with
polynomial dependence on $n$: since
$\ell\leq3(1+\ln n)$ and $(1+\ln n)^k\leq e\,k!\,n$,
we have $\ell^k\leq3^k e\,k!\,n$.
The second inequality follows from
$x^k\leq k!e^x$ for $x\in\mathbb R_{\geq0}$, applied to
$x=1+\ln n$.

The final output CNF is $\bigwedge_{i \in [1,m]} \Psi_i$. There are at most $\log_2 n$
nontrivial cyclic factors, and each clause contains at most $kn$
distinct variables. This proves the asserted size and construction
bounds and \eqref{eq:small-abelian-cnf}. 
\end{proof}

\subsection{An affine decomposition of the permutation action}
\label{sec:affine}

Let $n\in\N$ with $n\geq1$ and consider groups $A\unlhd G\leq\Sym(n)$ with $A$ abelian.
Let $O$ be the set of $A$-orbits on $[1,n]$.
Normality of $A$ implies that every $g\in G$ maps each $\Delta\in O$
onto another $A$-orbit, denoted by $\Delta^g\in O$.
For each $\Delta\in O$, let $A_\Delta=\{a\rest_\Delta:a\in A\}$ be the
permutation group induced by $A$ on $\Delta$.
Since $A$ is abelian, an element of $A$ fixing one point of $\Delta$
fixes every point of $\Delta$. Hence $A_\Delta$ acts regularly on $\Delta$
and $|A_\Delta|=|\Delta|$.
The disjoint union $\Sigma=\bigsqcup_{\Delta\in O}A_\Delta$ therefore has
$n$ elements.

Easdown and Praeger \cite[Proposition~1.3]{EasdownPraeger1988} show that,
if $A\neq1$ then there is an abelian normal subgroup $K\unlhd G$ containing
$A$ such that $G/K$ has a faithful permutation representation of degree
strictly smaller than $n$.
Seress \cite[Theorem~6.3.1]{Seress2003} describes an explicit action of
$G$ on $\Sigma$: for $g\in G$, $\Delta\in O$, and $a\in A$, the element
$a\rest_\Delta\in A_\Delta$ is sent to
$(g^{-1}ag)\rest_{\Delta^g}\in A_{\Delta^g}$.
One can choose the above subgroup $K$ as the kernel of this action.
Then $K$ is normal and abelian, and the action of $G$ on $\Sigma$ can
be constructed in polynomial time from generators for $G$ and $A$.
A related quotient construction is implemented in Magma as
\texttt{AbelianNormalQuotient} \cite{MagmaAbelianNormalQuotient}.

For each $\Delta\in O$, choose $o_{\scriptscriptstyle\Delta}\in\Delta$.
The bijection $A_\Delta\to\Delta$ given by $a\rest_\Delta\mapsto
o_{\!\scriptscriptstyle\Delta}^a$ for $a\in A$ identifies Seress's action with the
homomorphism $\lambda:G\to\Sym(n)$ from \cref{lem:affine}  below, as \eqref{eq:linear-part-conjugation} will show.
We give an additive formulation of this construction and include the proofs
to make the affine product formula \eqref{eq:affine-product} below explicit.
The operation
\begin{equation}
 o_{\!\scriptscriptstyle \Delta}^a + o_{\!\scriptscriptstyle \Delta}^b=o_{\!\scriptscriptstyle \Delta}^{ab}
 \qquad(a,b\in A)
 \label{eq:orbit-addition}
\end{equation}
defines an abelian group on $\Delta$ with identity $o_{\!\scriptscriptstyle \Delta}$; it is isomorphic to $A/\operatorname{Stab}_A(o_{\!\scriptscriptstyle \Delta})$,
where $\operatorname{Stab}_A(o_{\!\scriptscriptstyle \Delta}) = \{ a \in A : o_{\!\scriptscriptstyle \Delta}^a = o_{\!\scriptscriptstyle \Delta}\}$.
We call this group an \emph{orbit group}.
An orbit search yields for every $\alpha \in \Delta$ an $a_\alpha\in A$ with
$o_{\!\scriptscriptstyle \Delta}^{a_\alpha}=\alpha$. The permutations $a_\alpha$
allow us to compute the addition on $\Delta$ directly.
Hence, given generators for $G$ and $A$, for every orbit $\Delta \in O$ one can compute the addition table
of the orbit group $(\Delta, +)$ in polynomial time.

Define a permutation $\lambda(g) \in \Sym(n)$ for every $g \in G$ by
\begin{equation}
 \alpha^{\lambda(g)}=\alpha^g - o_{\!\scriptscriptstyle \Delta}^g \in \Delta^g
 \qquad(\Delta\in O,\ \alpha\in\Delta).
 \label{eq:linear-part}
\end{equation}
Note that $-$ refers to the abelian group on $\Delta^g$, since $\alpha^g, o_{\!\scriptscriptstyle \Delta}^g \in \Delta^g$.
We will in particular use the immediate consequence $\alpha^g = \alpha^{\lambda(g)}+ o_{\!\scriptscriptstyle \Delta}^g$,
which can be viewed as a decomposition of the action of $G$.
For any given $g \in G$ the permutation $\lambda(g)$ can be computed in polynomial time.

\begin{lemma}
\label{lem:affine}
The following hold for all $g \in G$, $\Delta \in O$ and $a \in A$:
\begin{itemize}
\item  $\left(o_{\!\scriptscriptstyle \Delta}^{a}\right)^{\lambda(g)}=(o_{\!\scriptscriptstyle \Delta^g})^{g^{-1}ag}$, in particular, $o_{\!\scriptscriptstyle \Delta}^{\lambda(g)}=o_{\!\scriptscriptstyle \Delta^g}$.
\item The mapping $\lambda:G\to\Sym(n)$ is a homomorphism
with $A\leq\ker\lambda$. 
\item For every $g \in G$, $\lambda(g) : (\Delta, +) \to (\Delta^g, +)$ is a group isomorphism.
\end{itemize}
\end{lemma}

\begin{proof}
Consider a $\Delta\in O$ and $g\in G$.
Choose $c\in A$ with
$o_{\!\scriptscriptstyle \Delta}^g=(o_{\scriptscriptstyle \Delta^g})^{c}$,
which is possible since both points belong to $\Delta^g$.
For $a\in A$, \eqref{eq:linear-part} gives
\begin{align}
 \left(o_{\!\scriptscriptstyle \Delta}^{a}\right)^{\lambda(g)}
 \stackrel{\text{\eqref{eq:linear-part}}}{=}
 \left(o_{\!\scriptscriptstyle \Delta}^{g}\right)^{g^{-1}ag}
       - o_{\!\scriptscriptstyle \Delta}^{g}
 =(o_{\!\scriptscriptstyle \Delta^g})^{cg^{-1}ag}
       - (o_{\scriptscriptstyle \Delta^g})^{c}
  \stackrel{\text{\eqref{eq:orbit-addition}}}{=} (o_{\!\scriptscriptstyle \Delta^g})^{cg^{-1}agc^{-1}}
  =(o_{\!\scriptscriptstyle \Delta^g})^{g^{-1}ag}.
 \label{eq:linear-part-conjugation}
\end{align}
Here $g^{-1}ag\in A$ by normality, and the last equality uses
that $A$ is abelian. 
 Taking $a=1\in A$ yields
 $o_{\!\scriptscriptstyle \Delta}^{\lambda(g)}=o_{\!\scriptscriptstyle \Delta^g}$.

Next, consider $\alpha,y\in\Delta$ and choose $a,b\in A$
with $\alpha=o_{\!\scriptscriptstyle \Delta}^{a}$ and
$y=o_{\!\scriptscriptstyle \Delta}^{b}$. By
\eqref{eq:linear-part-conjugation}, we obtain
\begin{align*}
 \left(\alpha+y\right)^{\lambda(g)}
 =\left(o_{\!\scriptscriptstyle \Delta}^{ab}\right)^{\lambda(g)}
  =o_{\!\scriptscriptstyle \Delta^g}^{(g^{-1}ag)(g^{-1}bg)}
 =o_{\!\scriptscriptstyle \Delta^g}^{g^{-1}ag}
       + o_{\!\scriptscriptstyle \Delta^g}^{g^{-1}bg}
  =\alpha^{\lambda(g)}+y^{\lambda(g)}.
\end{align*}
Thus the restriction of $\lambda(g)$ to $\Delta$ is a group
homomorphism from $\Delta$ to $\Delta^g$.
Since it is injective and $|\Delta|=|\Delta^g|$, it is an isomorphism.

We finally show that $\lambda : G \to \Sym(n)$ is a homomorphism. Let $g,h\in G$,
$\Delta\in O$, and $\alpha\in\Delta$. Choose $a\in A$ with
$\alpha=o_{\!\scriptscriptstyle \Delta}^{a}$. Since $g^{-1}ag\in A$ and
$(\Delta^g)^h=\Delta^{gh}$, repeated application of
\eqref{eq:linear-part-conjugation} gives
\begin{align*}
 \alpha^{\lambda(g)\lambda(h)}
 =\left((o_{\!\scriptscriptstyle \Delta^g})^{g^{-1}ag}\right)^{\lambda(h)}
  =(o_{\!\scriptscriptstyle \Delta^{gh}})^{h^{-1}g^{-1}agh}
 =(o_{\!\scriptscriptstyle \Delta^{gh}})^{(gh)^{-1}a(gh)}
  =\alpha^{\lambda(gh)}.
\end{align*}
Since the $A$-orbits partition $[1,n]$, this proves
$\lambda(gh)=\lambda(g)\lambda(h)$ for all $g,h\in G$.
For every $a\in A$, \eqref{eq:linear-part-conjugation} shows that $\lambda(a)$ is the identity permutation. 
\end{proof}

\begin{lemma}
\label{lem:affine-product}
Let $\ell\in\N$, $ g_1,\ldots, g_\ell\in G$, and $g= g_1\cdots  g_\ell\in G$.
For $i\in[1,\ell]$, put $p_i= g_1\cdots  g_{i-1}\in G$ and
$s_i= g_{i+1}\cdots  g_\ell\in G$. Then, for
$\Delta\in O$ and $\alpha\in\Delta$,
\begin{equation}
 \alpha^g=\alpha^{\lambda(g)}
       + \sum_{i=1}^{\ell}
       \left(\left(o_{\!\scriptscriptstyle \Delta}^{\lambda(p_i)}\right)^{ g_i}
       \right)^{\lambda(s_i)}.
 \label{eq:affine-product}
\end{equation}
Here, summation is done in the abelian group $(\Delta^g,+)$.
\end{lemma}

\begin{proof}
We prove \eqref{eq:affine-product} by induction on $\ell$. Fix
$\Delta\in O$ and $\alpha\in\Delta$. The case $\ell=0$ is clear. 
Now let $\ell\geq1$ and assume the formula holds for $\ell-1$.
Put $h= g_1\cdots  g_{\ell-1}=p_\ell\in G$. For
$i\in[1,\ell-1]$, put $t_i= g_{i+1}\cdots  g_{\ell-1}\in G$,
so that $s_i=t_i g_\ell$. The prefixes $p_i$ are unchanged.
With \cref{lem:affine}, we obtain the following (with all summations in $\Delta^g$):
\begin{eqnarray*}
 \alpha^g
 &\stackrel{\text{\eqref{eq:linear-part}}}{=}&\left(\alpha^{h}\right)^{\lambda( g_\ell)}
       + (o_{\!\scriptscriptstyle \Delta^{h}})^{ g_\ell}\\
 &\stackrel{\text{ind.~hyp.}}{=}& \alpha^{\lambda(h)\lambda( g_\ell)}
       + \sum_{i=1}^{\ell-1}
       \left(\left(o_{\!\scriptscriptstyle \Delta}^{\lambda(p_i)}\right)^{ g_i}
       \right)^{\lambda(t_i)\lambda( g_\ell)}
       + (o_{\!\scriptscriptstyle \Delta^{h}})^{ g_\ell}\\
 &=&\alpha^{\lambda(g)}
       + \sum_{i=1}^{\ell-1}
       \left(\left(o_{\!\scriptscriptstyle \Delta}^{\lambda(p_i)}\right)^{ g_i}
       \right)^{\lambda(s_i)}
       + \left(o_{\!\scriptscriptstyle \Delta}^{\lambda(h)}\right)^{ g_\ell}\\
 &=&\alpha^{\lambda(g)}
       + \sum_{i=1}^{\ell}
       \left(\left(o_{\!\scriptscriptstyle \Delta}^{\lambda(p_i)}\right)^{ g_i}
       \right)^{\lambda(s_i)}.
\end{eqnarray*}
This proves \eqref{eq:affine-product}.
\end{proof}

\subsection{Constructing circuits for product evaluation in solvable groups}
\label{sec:solvable-circuits}

In this section we fix a constant $d$. 
Moreover, for the following construction we consider integers 
$k,m,n\ge 1$, a permutation group $G \leq \Sym(n)$ with $G^{(d)}=1$
and an $m$-tuple $\overline{g} = (g_1, \ldots, g_m) \in G^m$.
We define the set of Boolean variables
\[ X_{k,m}=\{x_{i,j}:i\in[1,k],\ j\in[1,m]\}.
\] 
With $\mathcal I_k \subseteq 2^{[1,k]}$ we denote the set consisting of $\emptyset$ and all intervals
$[\ell,r]$ with $\ell,r\in[1,k]$ and $\ell\leq r$. A truth assignment $\phi : X_{k,m} \to \{0,1\}$
such that for every $i \in[1,k]$ there is a unique $j_i \in [1,m]$ with $\phi(x_{i,j}) = \mathbf{1}[j=j_i]$
is called \emph{admissible}. For such a $\phi$ and an interval $I \in \mathcal I_k$ we define the group element
\[
\overline{g}[I,\phi] = \begin{cases}
g_{j_\ell} g_{j_{\ell+1}}\cdots g_{j_r}\in G & \text{ if $I = [\ell,r]$ with $\ell,r\in[1,k]$ and $\ell\leq r$.}\\
1 & \text{ if } I = \emptyset .
\end{cases}
\]

\begin{lemma}
\label{lem:evaluation}
Given $n,m,k\ge 1$, $G$ and $\overline{g}=(g_1,\ldots,g_m)\in G^m$ as above, 
one can construct a Boolean circuit $C$ with input variables from $X_{k,m}$ and distinguished gates
 $B[I,\alpha,\beta]$ for all $I\in\mathcal I_k$ and $\alpha,\beta\in[1,n]$
 such that:\footnote{The circuit $C$ depends on $k,m,n, G$ and the tuple $(g_1,\ldots,g_m)$; this dependence is suppressed in the notation for the circuit.}
\begin{itemize}
\item For every admissible assignment $\phi : X_{k,m} \to \{0,1\}$,
$\phi(B[I,\alpha,\beta])=1$ if and only if $\alpha^{\overline{g}[I,\phi]}=\beta$.
\item  $C$ has weft at most $2d$ and depth at most $4d$. Every gate $B[I,\alpha,\beta]$ is either an input gate, a constant gate or an and-gate.
\item The total size and construction time for the circuit $C$ are $f_d(k)\cdot N^{c_d}$
for a computable function $f_d:\N\to\N$ and a constant
$c_d\in\N$, where $N$ is the total input length.
\end{itemize}
\end{lemma}

\begin{proof}
We construct the circuit by induction on $t\in[0,d]$ for permutation groups
satisfying $G^{(t)}=1$. At level $t$, the weft bound is $2t$ and
the depth bound is $4t$. 

For $t=0$ we have $G=1$, so every gate $B[I,\alpha,\beta]$ is a constant 
gate carrying the truth value $\mathbf{1}[\alpha=\beta]$. Now let $t\in[1,d]$ and suppose the assertion
holds for $t-1$. 

Compute first the abelian group $A=G^{(t-1)} \unlhd G$ and apply the construction from 
\cref{sec:affine}. We obtain the homomorphism $\lambda : G \to \Sym(n)$ with 
$A \leq \ker\lambda$ and the abelian group structures on the orbits of $A$.
Let $Q = G/\ker\lambda = \lambda(G)\leq\Sym(n)$. Since $A \leq \ker\lambda$, 
we have $Q^{(t-1)}=1$. Apply the induction hypothesis to $Q$ and the
tuple $(\lambda(g_1),\ldots,\lambda(g_m))$, with the same $n,m,k$.
Write $C'$ for the resulting circuit and $B'[I,\alpha,\beta]$ for its distinguished gates.
We will extend $C'$ to the circuit $C$ for $G$.

All gates $B[\emptyset,\alpha,\beta]$ in $C$ will be 
constant gates carrying the truth value $\mathbf{1}[\alpha=\beta]$.
Now assume that $I=[\ell,r]\in\mathcal I_k\setminus\{\emptyset\}$ and
$\alpha,\beta\in[1,n]$. For $i\in I$, define
$I_{<i}=[\ell,i-1]\in\mathcal I_k$ and
$I_{>i}=[i+1,r]\in\mathcal I_k$. Let $o\in[1,n]$ be the fixed representative 
 of the $A$-orbit containing $\alpha$. 
 For an admissible assignment $\phi$ such that $\phi(x_{i,j_i})=1$ (with $i \in [1,k]$, $j_i \in [1,m]$)
 we define points $\beta_i\in[1,n]$ for $i\in\{0\}\cup I$ by
\begin{equation}
 \beta_0=\alpha^{\lambda(\overline{g}[I,\phi])} \quad\text{and}\quad
 \beta_i=\left(o^{\lambda(\overline{g}[I_{<i},\phi])g_{j_i}}\right)^{\lambda(\overline{g}[I_{>i},\phi])}  \quad\text{for } i\in I.
 \label{eq:circuit-affine-summands}
\end{equation}
By \cref{lem:affine-product}, these points all belong to
the $A$-orbit containing $\alpha^{\overline{g}[I,\phi]}$, and their sum in that orbit
group is $\alpha^{\overline{g}[I,\phi]}$.

Let $\Delta$ be the $A$-orbit containing $\beta$. Recall that $\Delta$ carries an abelian group structure.
Let $s=|I|+1\in[1,k+1]$. We use \cref{lem:abelian-cnf} to
construct the CNF $\Phi[\Delta,s,\beta]$.
Recall that its propositional variables are $y_{i,a}$, where $i\in[1,s]$ and $a\in \Delta$.
We rename these variables to $y_{i,a}$ for $i\in \{0\}\cup I$ and $a\in \Delta$.
We then construct the positive CNF 
\begin{equation} \label{eq:CNF-for-circuit}
\Phi[\Delta,s,\beta] \wedge \bigwedge_{i\in\{0\}\cup I} \left(\bigvee_{a\in \Delta}y_{i,a}\right)
\end{equation}
and view it as an unbounded fan-in circuit of depth two.
The output gate of this circuit is $B[I,\alpha,\beta]$ and its input gates are labelled with the variables $y_{i,a}$ ($i\in\{0\}\cup I$, $a\in \Delta$).
Let us call the $y_{i,a}$-labelled input gate $B_{i,a}$. We then connect this CNF-circuit with the circuit $C'$ by 
connecting the gates $B_{i,a}$ as follows:
\begin{eqnarray}
 B_{0,a} &=&B'[I,\alpha,a],   \label{eq:circuit-summand-substitution1} \\
 B_{i,a}&=&\bigvee_{\substack{j\in[1,m]\\b\in[1,n]}}
 \bigl(x_{i,j}\wedge B'[I_{<i},o,b]
       \wedge B'[I_{>i}, b^{g_j}, a]\bigr)
       \ \text{ for } i\in I.
 \label{eq:circuit-summand-substitution}
\end{eqnarray}
Here, successive or-gates of unbounded fan-in are merged into a single or-gate. This happens when \eqref{eq:circuit-summand-substitution} replaces the variable 
$y_{i,a}$ in \eqref{eq:CNF-for-circuit}.

By the induction hypothesis, $B_{i,a}$ is true exactly when
$a=\beta_i$ in \eqref{eq:circuit-affine-summands}.
Hence the added
clauses in \eqref{eq:CNF-for-circuit} hold exactly when every $\beta_i$ belongs to $\Delta$. Under this
condition, for each position $i\in\{0\}\cup I$, exactly one
$a\in \Delta$ makes $B_{i,a}$ true, so \eqref{eq:small-abelian-cnf}
applies. Hence, the gate $B[I,\alpha,\beta]$ is true
exactly when all summands lie in $\Delta$ and sum to $\beta$, which is
equivalent to $\alpha^{\overline{g}[I,\phi]}=\beta$. 

We now verify the bounds on the weft and depth. 
For $t=0$, the depth and weft of the circuit are zero.
If $t>0$ then the circuit $C$ starts with a (possibly large) and-gate
followed by a (possibly large) or-gate, both coming from the positive CNF in \eqref{eq:CNF-for-circuit}.
Then the circuit either directly enters the circuit $C'$ (of depth at most $4(t-1)$ and weft at most $2(t-1)$) by \eqref{eq:circuit-summand-substitution1}
or it reaches the input variable $x_{i,j}$ in \eqref{eq:circuit-summand-substitution} or it traverses at most two binary and-gates before
entering $C'$.  This yields a total depth of $4+4(t-1) = 4t$ and weft $2+2(t-1)=2t$.

The size of the circuit can be estimated as follows (we count circuit wires, which is more appropriate than counting gates for unbounded fan-in circuits): 
Going from the circuit $C'$ to the circuit $C$ adds
$(k(k+1)/2+1)n^2 = \Theta(k^2 n^2)$ new gates $B[I,\alpha,\beta]$.  For each of these gates we 
construct the CNF $\Phi[\Delta,s,\beta]$ of size $f(k) \cdot n^c$ (for a computable
function $f:\N\to\N$ and a constant $c$) from \cref{lem:abelian-cnf}.
The additional clauses in \eqref{eq:CNF-for-circuit} yield a CNF of size $f(k) \cdot n^c + \mathcal{O}(k n)$.
After adding the circuit wires from \eqref{eq:circuit-summand-substitution} and merging successive or-gates
we obtain $\mathcal{O}(m n(f(k) n^c + kn))$ new wires in the circuit. This has to be done for each of the $\Theta(k^2 n^2)$
many new gates $B[I,\alpha,\beta]$. In total, we add $\mathcal{O}(m f(k) n^c)$ new wires in the circuit for a new
computable function $f$ and a new constant $c$.

This is repeated $d$ times, where $d$ is a constant in our setting.
In total, the circuit size is bounded by $f_d(k)N^{c_d}$ for
a computable function $f_d:\N\to\N$ and a constant
$c_d\in\N$ depending only on $d$.
The same bound holds for the construction time.
The addition tables of the orbit groups
and cyclic coordinates are computable in polynomial time.
\end{proof}

\begin{theorem}
\label{thm:solvable}
For every fixed $d\ge 1$, the problem
$\pSSS[\SolPG{d}]$ belongs to
$\W[2d]$.
\end{theorem}

\begin{proof}
Fix an instance of $\pSSS[\SolPG{d}]$, consisting of 
$g_1,\ldots,g_m,h\in \Sym(n)$ with $G := \langle g_1,\ldots,g_m \rangle$ and $G^{(d)}=1$,
 and the parameter $k\in [1,m]$. W.l.o.g.~we can assume that $h \in G$.
We construct the circuit $C$ from \cref{lem:evaluation} 
with distinguished gates $B[I,\alpha,\beta]$ and input variables $x_{i,j}$ ($i \in [1,k]$, $j \in [1,m]$)
and extend $C$ by an output gate computing
\begin{equation*}
 \bigwedge_{i=1}^k\bigvee_{j=1}^m x_{i,j} \ \wedge
 \bigwedge_{\substack{i\in[1,k-1]\\j,j'\in[1,m],\ j\geq j'}} \!\!\!\!\!\!\!\!\!
      (\neg x_{i,j}\vee\neg x_{i+1,j'})
\ \wedge \ \bigwedge_{\alpha\in[1,n]} \!\!\! B[[1,k],\alpha,\alpha^h],
 \label{eq:solvable-weighted-circuit}
\end{equation*}
where the last conjunction is merged with a possible top and-gate in $B[[1,k],\alpha,\alpha^h]$.

The first conjunction,
together with the restriction to truth assignments of weight exactly $k$, forces 
for every $i \in [1,k]$ exactly one variable $x_{i,j}$ ($j \in [1,m]$) to be true. The second conjunction forces these
indices to be strictly increasing. The final conjunction then
holds exactly when their ordered product is $h$.

Since every gate $B[[1,k],\alpha,\alpha^h]$ is an and-gate, a constant gate, or an input gate in the circuit
$C$, the circuit resulting from the above construction still has depth $4d$ and weft $2d$ (if $d \geq 1$).
Our $\pSSS[\SolPG{d}]$ instance is positive if and only if the circuit has a satisfying assignment of weight $k$.
\end{proof}
\cref{thm:solvable} gives the upper bounds $\W[2]$ for abelian
permutation groups and $\W[4]$ for metabelian permutation groups.
\cref{pre:lem:replication,pre:lem:padding} show that the problems 
$\pSSS_\le[\SolPG{d}]$, $\pF_\le[\SolPG{d}]$ and $\pF[\SolPG{d}]$ belong to
$\W[2d]$ as well.

\section{Nilpotent permutation groups of fixed class}
\label{sec:nilpotent}

In this section, $d \ge 1$ will again be a fixed constant. 
We consider the restriction of the problem $\pSSS[\Sym(*)]$, where the input permutations $g_1, \ldots, g_m$ generate a
$d$-step nilpotent permutation group. Recall from \cref{sec:PGs} that this property can be checked in polynomial time.
We denote this restriction of $\pSSS[\Sym(*)]$ with $\pSSS[\NilPG{d}]$; the problems $\pSSS_\le[\NilPG{d}]$, $\pF[\NilPG{d}]$ and
$\pF_\le[\NilPG{d}]$ are defined similarly. Note that $\pSSS[\NilPG{1}]$ is the same problem as $\pSSS[\SolPG{1}]$, which is the restriction of 
 $\pSSS[\Sym(*)]$ to abelian permutation groups. 
 
 The main result of this section states that $\pSSS[\NilPG{d}]$ belongs to $\W[2]$ for every $d \geq 1$. We start with some well-established
 concepts for free nilpotent groups.

\subsection{Free nilpotent groups, Hall coordinates and the Magnus embedding} \label{sec:Hall-Magnus}

The main sources for this section are \cite{Hall1969,Reutenauer1993}.
The \emph{free $d$-step nilpotent
group} of rank $r$ with generators $a_1,\ldots,a_r$ is the quotient
$\mathcal{N}_{d,r}=\mathcal{F}_r/\gamma_{d+1}(\mathcal{F}_r)$, where $\mathcal{F}_r = \mathcal{F}_r(a_1,\ldots,a_r)$ is the free group 
generated by $a_1,\ldots,a_r$. In our later setting, $r$ will be part of the input but $d$ will always be fixed and not part of the input.
Let us define the alphabet $\Sigma_r = \{ a_1, a_1^{-1},\ldots,a_r,a_r^{-1} \}$.
We define the \emph{basic commutators} and their order inductively following \cite[Section~5]{Hall1969} and \cite[Section~4.1, (4.1.1)--(4.1.3)]{Reutenauer1993}.
For this we first define  \emph{basic expressions} as formal expressions built up from the 
symbols $a_1,\ldots,a_r$ and the binary operator $[\cdot, \cdot]$ and thereby
define a linear order on them.
The generators $a_1<\cdots<a_r$ are basic expressions of weight~$1$.
For $n\in [2,d]$ and previously constructed basic expressions $u, v$ of weights $\ell, m$ with $n = \ell+m$,
the expression $[u, v]$ is basic of weight $n$ if $u > v$ and the following condition holds:
if $u = [u_1,u_2]$  then $u_2\leq v$.
We order basic expressions by increasing weight; those of the same weight greater than $1$ are ordered arbitrarily.
Note that basic expressions can be viewed as ordered binary trees whose leaves are labelled with symbols from $a_1,\ldots,a_r$; they
are called Hall trees in \cite{Reutenauer1993}.
Note that Reutenauer's order on Hall trees is the reverse of our order on basic expressions.

 We obtain basic commutators of weight $s \in [1,d]$ by interpreting basic expressions of weight $s$ in $\mathcal{N}_{d,r}$.
 Distinct basic expressions of weight at most $d$ represent distinct basic commutators in $\mathcal{N}_{d,r}$ by Hall’s normal form theorem recalled below. 
 Moreover, basic commutators inherit the linear order from basic expressions. Note that every basic commutator $a_i$  belongs
 to $\gamma_{\omega_i}(\mathcal{N}_{d,r})$; this follows from the general inclusion $[\gamma_i(G),\gamma_j(G)]\subseteq\gamma_{i+j}(G)$.

Extend the ordered list of generators $a_1,\ldots,a_r$ to the list
$a_1,\ldots,a_t\in \mathcal{N}_{d,r}$ of all basic commutators of weights at most
$d$, listed in order, and let $\omega_i$ be the weight of $a_i$.
One has $t = \mathcal{O}(r^d)$ by Witt's formula \cite[Theorem~5.7]{Hall1969}.
The importance of basic commutators comes from the following result of Hall;
see e.g. \cite[Theorem 6.8]{Reutenauer1993}.

\begin{theorem}[Hall’s normal form theorem] \label{thm:Hall-nf}
Every $g\in \mathcal{N}_{d,r}$ can be uniquely written as
$g=a_1^{y_1}\cdots a_t^{y_t}$ with
$y_1,\ldots,y_t\in\Z$.
\end{theorem}
We call the integers $y_i$ the \emph{Hall coordinates} of $g$ 
and write $\Hall_i(g) = y_i$. The process of bringing a word over the alphabet $\Sigma_r$
into the above Hall normal form is known as commutator collection.
An element $g \in \mathcal{N}_{d,r}$ belongs to $\gamma_s(\mathcal{N}_{d,r})$ if 
and only if $\Hall_i(g) = 0$ for all $i \in [1,t]$ with $\omega_i < s$ \cite[Theorem~6.8 and Corollary~6.15]{Reutenauer1993}.

The second important concept for nilpotent groups that we use is the Magnus embedding.
For non-commuting variables $\chi_1, \ldots, \chi_r$ we consider the 
free algebra $\Q\langle \chi_1,\ldots, \chi_r\rangle$ of all non-commutative polynomials over the field $\Q$.
Monomials are then simply words over the alphabet $\{\chi_1, \ldots, \chi_r\}$.
With $\mathfrak m$ we denote the two-sided ideal generated by $\chi_1,\ldots,\chi_r$.
Thus $\mathfrak m^{d+1}$ is the ideal generated by all words in $\{\chi_1,\ldots,\chi_r\}^{d+1}$.
The truncated free associative algebra $\mathcal{A}_{d,r}$ is defined as the quotient algebra
\[
 \mathcal{A}_{d,r} = \Q\langle \chi_1,\ldots, \chi_r\rangle/\mathfrak{m}^{d+1}.
 \]
 Every element of $\mathcal{A}_{d,r}$ can be uniquely written as a $\Q$-linear combination of words
 from $\{ \chi_1, \ldots, \chi_r \}^{\le d}$.
 Hence, as a vector space over $\Q$, $\mathcal{A}_{d,r}$ 
has dimension $\Theta(r^d)$. For a word $w\in\{\chi_1,\ldots,\chi_r\}^*$ and $P \in \mathcal{A}_{d,r}$
we write  $[w]P \in \Q$ for the coefficient of $w$ in $P$.

The \emph{Magnus embedding} $\mu_{d,r} : \mathcal{N}_{d,r} \to \mathcal{A}_{d,r}$
 is defined by $\mu_{d,r}(a_j)=1+\chi_j$ for $j\in[1,r]$. It maps $\mathcal{N}_{d,r}$ injectively into the group of units
of $\mathcal{A}_{d,r}$ \cite[Theorem~5.6]{Hall1969}, which consists of all non-commutative polynomials whose constant term is not zero.
In our application the constant term will always be $0$ or $1$ and all other coefficients will be from $\Z$.
For  $U\in \mathcal{A}_{d,r}$ with zero constant term, we have
$(1+U)^{-1}=\sum_{\ell=0}^d(-U)^\ell$, since $U^{d+1}=0$. 

By induction on $s$ it follows easily that, for every $g \in \gamma_s(\mathcal{N}_{d,r})$, every monomial of $\mu_{d,r}(g)-1 \in \mathcal{A}_{d,r}$ has 
length at least $s$ \cite[Lemma~6.10(iii)]{Reutenauer1993}. In particular, for every basic commutator $a_i$ ($i \in [1,t]$) we have
$\mu_{d,r}(a_i) = 1 + A_i$, where all monomials of $A_i$ have length at least $\omega_i$ \cite[Lemma~6.10(ii)]{Reutenauer1993}.

\begin{lemma} \label{claim:size-coefficients}
Consider a product $U = \prod_{i=1}^k (1+U_i)$ with $k \geq 1$, where every $U_i \in \mathcal{A}_{d,r}$ has constant term $0$ and every coefficient in $U_i$ is an integer
of bit length at most $\sigma$. Then, every coefficient of $U$ is an integer of bit length $\mathcal{O}(\log k + \sigma)$.
\end{lemma}

\begin{proof}
Assume that all coefficients of the $U_i$ have absolute value at most $\alpha \ge 1$.
Consider a monomial $w$ in $U$; it is a word of length at most $d$. Hence, in the expansion of the product $\prod_{i=1}^k (1+U_i)$, the monomial can be obtained
in at most $k^{d+1}$ ways. Hence, the absolute value of the coefficient of $w$ is bounded by $k^{d+1} \alpha^d$ and its bit length is
bounded by $\mathcal{O}(\log k + \log \alpha) = \mathcal{O}(\log k + \sigma)$
(recall that $d$ is a constant).
\end{proof}

\subsection{The main lemma}

The following lemma is the main technical result of this section. We assume that 
polynomials in $\Z[y_1,\ldots,y_m]$ are given by their sparse representation, which stores
all monomials together with the corresponding coefficient in binary notation. The value $d$ is considered as a fixed constant
throughout this section. All constants hidden in the $\mathcal{O}$-notation only depend on $d$.

\begin{lemma}
\label{lem:nilpotent-polynomials}
Given $m,n\ge 1$ in unary notation,
a $d$-step nilpotent group $G\leq\Sym(n)$ and elements
$g_1,\ldots,g_m,h\in G$, one can compute in polynomial time
the binary encoding of an integer $q \ge 1$ together with its prime factorization and
polynomials $P_1,\ldots,P_t\in\Z[y_1,\ldots,y_m]$ of degree at most $d$
such that $t=\mathcal{O}(m^d)$ and for every $\overline{\epsilon} = (\epsilon_1, \ldots, \epsilon_m)
\in\{0,1\}^m$,
\begin{equation}
 g_1^{\epsilon_1}\cdots g_m^{\epsilon_m}=h
 \quad\Longleftrightarrow\quad
\forall i\in[1,t] : P_i(\overline{\epsilon})\equiv0\pmod q.
 \label{eq:nilpotent-polynomials}
\end{equation}
If $q = p_1^{v_1} \cdots p_l^{v_l}$ is the prime factorization of $q$ then $l = \poly(n)$ and
$p_i^{v_i} = \poly(n)$ for all $i \in [1,l]$.
\end{lemma}
The existence of polynomial congruences of degree at most $d$ follows
from \cite[Theorem~9.1]{Podkorytov2026}. We give a construction that also satisfies
the additional bounds from \cref{lem:nilpotent-polynomials}.
We first prove several auxiliary lemmas.

For $r\geq1$, use the basic commutators $a_1,\ldots,a_t\in\mathcal{N}_{d,r}$
and their weights $\omega_i\in[1,d]$ from \cref{sec:Hall-Magnus}.
All products indexed by a subset $I \subseteq [1,t]$ are taken in increasing order on $I$, unless stated otherwise.
For $s\in[1,d]$, define the subintervals $I_s=\{i\in[1,t]:\omega_i=s\}$ and $I_{<s}=\{i\in[1,t]:\omega_i < s\}$
of $[1,t]$ and let
$r_s=|I_s|\in\N$. The intervals $I_{\le s}$ and $I_{\ge s}$ are defined analogously. Note that $r_1=r$.
 In the following, for a tuple $\overline{u} = (u_1, \ldots, u_t)$ (with the $u_i$ from any set)
we write $\overline{u}_s=(u_i)_{i\in I_s}$ and
$\overline{u}_{<s}=(u_i)_{i \in I_{<s}}$. The same notation is used for column vectors.

We start with showing that from an element $g \in \mathcal{N}_{d,r}$ that is given by an SLP over
the standard generators, one can compute the Hall coordinates of $g$ in polynomial time.

\begin{lemma}
\label{lem:slp-hall-coordinates}
Given $r\geq1$ in unary notation and an SLP of size
$\ell\geq1$ over $\Sigma_r$ representing
$g\in\mathcal{N}_{d,r}$, one can compute all Hall coordinates of $g$
in time $\poly(r,\ell)$.
\end{lemma}

\begin{proof}
We use the Magnus embedding $\mu_{d,r}$ introduced in \cref{sec:Hall-Magnus}.
In the following, we write $\mu$ for $\mu_{d,r}$ and $\mathcal{A}$ for
$\mathcal{A}_{d,r}$.  The list $a_1,\ldots,a_t$ can be computed in polynomial time in $r$,
if we allow the commutator bracket $[\cdot,\cdot]$ as a syntactic operator.
Every $a_i$ can be written as a word over $\Sigma_r$
of length depending only on $d$.
Replacing in this word every symbol $a_j$ ($j \in [1,r]$) by $1+\chi_j$ and every symbol
$a_j^{-1}$ by $1-\chi_j+\chi_j^2-\cdots+(-1)^d\chi_j^d$ and expanding the product
yields $\mu(a_i)=1+A_i$. Each $A_i\in\mathcal{A}$ has $\mathcal{O}(1)$
monomials of length at least $\omega_i$, and all its coefficients have
bit length $\mathcal{O}(1)$ by \cref{claim:size-coefficients}.

From the SLP for $g$ one can compute the polynomial $\mu(g)$ explicitly.
All monomials in $\mu(g)$ are words over $\{\chi_1,\ldots,\chi_r\}$ of length
at most $d$. Hence, for each SLP instruction we have to multiply two
polynomials with $\mathcal{O}(r^d)$ possible monomials.
The SLP can be expanded into a word over $\Sigma_r$ of length at most $2^\ell$.
Replacing every symbol $a_j^{\pm 1}$ ($j \in [1,r]$) by $\mu(a_j^{\pm 1})$ and applying
\cref{claim:size-coefficients} shows that all coefficients of $\mu(g)$ have bit length
$\mathcal{O}(\ell)$. The same bound holds for every intermediate SLP result,
so evaluating the SLP in $\mathcal{A}$ takes time $\poly(r,\ell)$.

We now describe how to recover the Hall coordinates of $g\in \mathcal{N}_{d,r}$
from $\mu(g)\in \mathcal{A}$.
For $i\in[1,t]$, let $B_i\in \mathcal{A}$ be the homogeneous component of degree
$\omega_i$ of $\mu(a_i)-1 = A_i$. Thus $B_i$ is obtained simply by retaining the monomials of length $\omega_i$
in the explicitly computed polynomial $\mu(a_i)-1$. 
By \cite[Theorem~4.9(i) and Lemma~6.10(ii)]{Reutenauer1993},  for every $s\in[1,d]$
the polynomials $B_i$ for $i\in I_s$ are linearly independent over $\Q$; in particular, they are pairwise distinct.

Let $z_i=\Hall_i(g)\in\Z$ for $i\in[1,t]$ so that
$g=a_1^{z_1}\cdots a_t^{z_t}$. We recover the coordinates $z_i$ from $\mu(g)$ successively in blocks for each
weight $s\in[1,d]$. At the beginning of stage $s$, we maintain a polynomial
$V_s\in \mathcal{A}$ satisfying
\begin{equation}
 V_s=\mu\bigg(\prod_{i \in I_{\ge s}} a_i^{z_i}\bigg).
 \label{eq:hall-recovery-residual}
\end{equation}
Initially, $V_1=\mu(g)$, so \eqref{eq:hall-recovery-residual} holds for $s=1$.
For $s\in[1,d]$, let $U_s\in \mathcal{A}$ be the homogeneous length-$s$ component of $V_s-1$. We have
\begin{equation}
 \sum_{i\in I_s}z_iB_i = U_s.
 \label{eq:hall-recovery-linear-system}
\end{equation}
Indeed, for $i\in I_s$, the polynomial $\mu(a_i^{z_i})$ has the form $1 + z_iB_i + W_i$ for $W_i\in\mathcal{A}$, where 
all monomials in $W_i$ have length larger than $s$. Multiplying these polynomials and taking the 
homogeneous length-$s$ part yields \eqref{eq:hall-recovery-linear-system}.

Equating coefficients of the $r^s$ monomials $w \in \{\chi_1,\ldots,\chi_r\}^s$
turns \eqref{eq:hall-recovery-linear-system} into a linear system for the
$r_s=|I_s| \leq r^s$ unknowns $z_i$ for $i\in I_s$. The right-hand side of the equation for $w\in \{\chi_1,\ldots,\chi_r\}^s$
is the coefficient of $w$ in $U_s$. The matrix $C_s\in\Z^{r^s\times r_s}$ of the linear system 
has in row $w\in\{\chi_1,\ldots,\chi_r\}^s$ and column $i\in I_s$
the entry $[w]B_i$. The matrix $C_s$ has full column rank by the linear independence
of the $B_i$. Moreover, the system has the integer solution $z_i = \Hall_i(g)$ for $i\in I_s$, 
which must be the unique solution over  $\Q$.

If $r_s=0$, there are no coordinates $z_i$ to recover at this stage and
we set $V_{s+1}=V_s$. Assume now that $r_s\geq1$ and that every coefficient
of $V_s$ has bit length at most $\alpha\in\N$ with $\alpha\geq1$.
Select $r_s$ linearly independent rows of $C_s$ to obtain an invertible matrix
$D_s\in\Z^{r_s\times r_s}$. The absolute value of every entry of $D_s$ is bounded by a constant $c \geq1$ depending only on $d$,
whereas the entries of the right-hand side vector are the coefficients of $U_s$ (and zero) and therefore 
have absolute value at most
$2^\alpha$. For $i\in I_s$, let $D_{s,i}\in\Z^{r_s\times r_s}$
be obtained by replacing in $D_s$ the $i$-th column by the  right-hand side vector.
Since $|\det D_s|\geq1$, Cramer's rule and Hadamard's inequality give
$|z_i|\leq|\det D_{s,i}|\leq2^\alpha r_s^{r_s/2}c^{r_s-1}$.
Thus every $z_i$ has bit length at most
$\alpha+\mathcal{O}(r_s\log r_s) \le \alpha+\mathcal{O}(r^s\log r)$.
The solution $(z_i)_{i\in I_s}$ can be computed in time $\poly(r,\alpha)$ by Gaussian elimination over $\Q$
\cite[Theorem~3.3]{Schrijver1986}.

Having recovered the Hall coordinates $z_i$ ($i \in I_s$), we compute the invertible polynomial
\[ V'_s=\mu\bigg(\prod_{i\in I_s}a_i^{-z_i}\bigg) \in\mathcal{A},\] 
where the product runs over the set $I_s$ in decreasing order, and set
$V_{s+1}=V'_sV_s$. This yields
$V_{s+1}=\mu(\prod_{i\in I_{\ge s+1}}a_i^{z_i})$, so
\eqref{eq:hall-recovery-residual} holds at the next stage.
After stage $d$, all Hall coordinates have been recovered and $V_{d+1}=1$.

It remains to bound the coefficient lengths during these updates.
For every $i\in I_s$, the finite binomial expansion
$\mu(a_i^{-z_i})=\sum_{j=0}^{\lfloor d/s\rfloor}\binom{-z_i}{j}A_i^j$
holds for every $z_i\in\Z$.
Since $d$ is fixed and each $A_i$ has constantly many monomials with
coefficients of constant bit length, this computes every factor
$\mu(a_i^{-z_i})$ with coefficient bit length
$\mathcal{O}(\alpha+r^d\log r)$.
The polynomial $V_{s+1}$ is the product of the polynomials $\mu(a_i^{-z_i})$ (in decreasing order over $I_s$)
and $V_s$. Each of these $r_s+1 \leq r^s+1$ polynomials has coefficients of bit length bounded by $\mathcal{O}(\alpha+r^d\log r)$.
By \cref{claim:size-coefficients}, all coefficients of
this product and its partial products have bit length
$\mathcal{O}(\alpha+r^d\log r)$. There are $\mathcal{O}(r^d)$ possible monomials and at most
$r^s+1 \leq r^d+1$ polynomials that have to be multiplied. Hence, the running time for computing $V_{s+1}$ from $V_s$ is
 $\poly(r,\alpha)$.

Recall that we start with $V_1 = \mu(g)$, whose coefficients have bit length $\mathcal{O}(\ell)$.
There are only $d = \mathcal{O}(1)$ stages, and in each stage the bit length increases from $\alpha$ to $\mathcal{O}(\alpha+r^d\log r)$.
Hence all intermediate coefficients and all
Hall coordinates have bit length $\mathcal{O}(\ell+r^d\log r)= \poly(r,\ell)$.
This proves the claimed polynomial running time.
\end{proof}
For the remaining lemmas, fix the input data of \cref{lem:nilpotent-polynomials}.
Polynomial time will always mean time $\poly(m,n)$.
Let $r=m+1\in\N$ and let $\mathcal{N}_{d,r}$ be the free $d$-step nilpotent group generated by $a_1,\ldots,a_r\in \mathcal{N}_{d,r}$.
Let $H=\langle g_1,\ldots,g_m,h\rangle\leq G$, and let
$\pi:\mathcal{N}_{d,r}\to H$ be the homomorphism with $\pi(a_i)=g_i$ for
$i\in[1,m]$ and $\pi(a_r)=h$. Let $K=\ker\pi\leq \mathcal{N}_{d,r}$.

We use the basic commutators $a_1, \ldots, a_t$ and index sets  $I_s$, $I_{<s}$, etc. defined above; in particular,
$t=\mathcal{O}(r^d)=\mathcal{O}(m^d)$.
For $i\in[1,t]$, let $N_i=\langle a_i,\ldots,a_t\rangle\leq\mathcal{N}_{d,r}$
and $H_i=\pi(N_i)\leq H$, and set $N_{t+1}=H_{t+1}=1$.
Let $e_i=[H_i:H_{i+1}]\in[1,n!]$ for $i\in[1,t]$.

From the corollary to \cite[Theorem~5.6]{Hall1969} one can easily derive that the groups $N_i$ are normal in $\mathcal{N}_{d,r}$, and that $N_i/N_{i+1}$ is infinite
cyclic, generated by $a_iN_{i+1}$.  We present a proof of these facts to keep the paper self-contained.

\begin{lemma} \label{lemma:N_i}
Every group $N_i$ is normal in $\mathcal{N}_{d,r}$, and $N_i/N_{i+1}$ is infinite
cyclic, generated by $a_iN_{i+1}$. For every $g \in N_i$ we have $\Hall_j(g)=0$ for $j \in [1,i-1]$.
\end{lemma}

\begin{proof}
Fix an $i\in[1,t]$, let $s=\omega_i\in[1,d]$, and let $p\in[i,t]$
be the largest index with $\omega_p=s$. Consider the canonical epimorphism
$\eta:\mathcal{N}_{d,r}\to\mathcal{N}_{s,r}$, and write
$\hat{a}_j=\eta(a_j)\in\mathcal{N}_{s,r}$ for $j\in[1,t]$.
Then $\hat{a}_1,\ldots,\hat{a}_p$ is the  ordered list of
basic commutators of weight at most $s$ in $\mathcal{N}_{s,r}$,
whereas $\hat{a}_j=1$ for $j>p$.

Let $g=a_1^{y_1}\cdots a_t^{y_t}\in\mathcal{N}_{d,r}$, where
$y_1,\ldots,y_t\in\Z$ are the Hall coordinates of $g$. Then
$\eta(g)=\hat{a}_1^{y_1}\cdots\hat{a}_p^{y_p}$.
By uniqueness of the Hall normal form in $\mathcal{N}_{s,r}$ (\cref{thm:Hall-nf}), we have $\eta(g)=1$ if and
only if $y_1=\cdots=y_p=0$. Since $\eta(a_j)=1$ for $j>p$,
this proves $\ker\eta=N_{p+1}$, with $N_{t+1}=1$ as usual.

The elements $\hat{a}_i,\ldots,\hat{a}_p$ belong to
$\gamma_s(\mathcal{N}_{s,r})$ and are therefore central in
$\mathcal{N}_{s,r}$, which yields $\langle \hat{a}_i,\ldots,\hat{a}_p \rangle \unlhd \mathcal{N}_{s,r}$.
Moreover, $\ker\eta=N_{p+1}\subseteq N_i$.
Consequently,
$N_i=\eta^{-1}(\langle\hat{a}_i,\ldots,\hat{a}_p\rangle)$
is normal in $\mathcal{N}_{d,r}$.

Now suppose that $g=a_1^{y_1}\cdots a_t^{y_t}\in N_i$. Since
$\hat{a}_i,\ldots,\hat{a}_p$ commute, there are
$z_i,\ldots,z_p\in\Z$ such that
$\eta(g)=\hat{a}_i^{z_i}\cdots\hat{a}_p^{z_p}$.
Comparing this expression with
$\eta(g)=\hat{a}_1^{y_1}\cdots\hat{a}_p^{y_p}$,
uniqueness of the Hall normal form gives $y_j=0$ for $j\in[1,i-1]$, i.e., $\Hall_j(g)=0$ for all $j\in[1,i-1]$.

Finally, $N_i=\langle a_i,N_{i+1}\rangle$, so the quotient
$N_i/N_{i+1}$ is cyclic and generated by $a_iN_{i+1}$.
If $z\in\Z$ satisfies $a_i^z\in N_{i+1}$ then $\Hall_i(a_i^z)=0$ by the conclusion from the previous paragraph. On the other
hand, $\Hall_i(a_i^z)=z$ by uniqueness of Hall normal form.
Hence $z=0$, and $N_i/N_{i+1}$ is infinite cyclic.
\end{proof}
The following lemma shows that one can compute SLPs for elements $b_i\in K\cap N_i$
for $i\in[1,t]$ that establish a Hall normal form for every element from the kernel $K$:

\begin{lemma}
\label{claim:unique-decomp-b_i}
One can compute in polynomial time SLPs over
$\Sigma_r$ for elements $b_i\in K\cap N_i$ $(i\in[1,t])$ such that $\Hall_i(b_i)=e_i$ and $\Hall_j(b_i)=0$ for $j\in[1,i-1]$.
Every element $g\in K$ can be uniquely written as
$g=b_1^{s_1}\cdots b_t^{s_t}$ for $s_1,\ldots,s_t\in\Z$.
\end{lemma}

\begin{proof}
We have $H_i = \langle \pi(a_i),\ldots, \pi(a_t)\rangle$ and the generators $\pi(a_j)$ can be computed
in polynomial time by evaluating the commutator expressions for the $a_j$ in $\Sym(n)$.
For $i\in[1,t]$, we then compute the index $e_i=[H_i:H_{i+1}]\in [1,n!]$. This is possible 
in polynomial time, since the order of a permutation group can be computed in polynomial time \cite[Section~3.1, item~(c)]{Seress2003}.
Then $\pi(a_i)^{e_i}\in H_{i+1}$ and using fast exponentiation, we can compute the permutation $\pi(a_i)^{e_i}$ in 
polynomial time. By \cref{theorem:slp} one can produce in polynomial time an SLP of size $\poly(n)$ over
$\{1,\pi(a_{i+1}),\ldots,\pi(a_t)\}$ for
$\pi(a_i)^{e_i}$ ($1$ is needed if $i=t$).
We lift this SLP to an SLP for 
an element $v_{i+1}\in N_{i+1}$ with $\pi(v_{i+1}) = \pi(a_i)^{e_i}$
by replacing in the SLP every $\pi(a_j)$ by $a_j$ ($j \in [i+1,t]$). We then define
\[ b_i=a_i^{e_i}v_{i+1}^{-1}\in K\cap N_i.\]
We have $b_i N_{i+1} = a_i^{e_i} N_{i+1}$, $\Hall_i(b_i) = e_i$ and $\Hall_j(b_i) = 0$ for $j \in [1,i-1]$.
Moreover, we obtain in polynomial time an SLP for $b_i$ over the generators $a_i^{\pm 1}, \ldots, a_t^{\pm 1}$ of $N_i$ and its size is $\poly(n)$.
This SLP can finally be transformed in polynomial time into an equivalent SLP over  $\Sigma_r$,
by replacing every $a_j^{\pm 1}$ ($j \in [i,t]$) by an equivalent word over $\Sigma_r$. Since 
$a_j$ is a basic commutator of weight at most $d$, the length of this word only depends on the constant $d$.
The identity can be represented by $a_1a_1^{-1}$ if necessary.
Hence, the resulting SLP for $b_i$ still has size $\poly(n)$.

For $i\in[1,t]$, consider an element $g \in K\cap N_i$. It can be uniquely written
as $g = a_i^{y_i} a_{i+1}^{y_{i+1}} \cdots a_t^{y_t}$ for integers $y_j$.
Since $g \in K$ we have $\pi(a_i)^{-y_i} \in H_{i+1}$, i.e., $\pi(a_i)^{-y_i} = 1$ in $H_i/H_{i+1}$,
which implies that $y_i$ is a multiple of $e_i$. Let $q_i = y_i/e_i$. With the normality of $N_{i+1}$  we obtain
$(b_i^{-q_i}g)N_{i+1}=a_i^{-e_iq_i+y_i}N_{i+1}=N_{i+1}$,
i.e., $b_i^{-q_i}g \in N_{i+1}$. 
Combined with $b_i^{-q_i}g \in K$, this proves $b_i^{-q_i}g\in K\cap N_{i+1}$.

Continuing this way with $b_i^{-q_i} g$ (and using $N_{t+1} = 1$) proves the existence of the factorization from the lemma.
To see uniqueness, assume that $b_1^{s_1}\cdots b_t^{s_t} = b_1^{u_1}\cdots b_t^{u_t}$
and let $i$ be the smallest index with $s_i \neq u_i$. After cancelling the common initial factors, we obtain in the quotient $N_i/N_{i+1}$ an identity $b_i^s = 1$ for $s=s_i-u_i\neq 0$, which yields
$a_i^{e_i s} = 1$ in $N_i/N_{i+1}$. Since $N_i/N_{i+1}$ is infinite cyclic by \cref{lemma:N_i} and $e_i \neq 0$, this yields a contradiction, which proves 
\cref{claim:unique-decomp-b_i}.
\end{proof}
Introduce commuting variables $x_1,\ldots,x_t$, assign to $x_i$ the weight $\omega_i \ge 1$, and set
$\overline{x}=(x_1,\ldots,x_t)$. For a polynomial $P \in \Q[\overline{x}]$ we define the weighted degree of $P$ 
as follows: For a monomial $x_1^{n_1} x_2^{n_2} \cdots x_t^{n_t}$ the weighted degree is $n_1 \omega_1 + n_2 \omega_2 + \cdots + n_t \omega_t$
and the weighted degree of $P$ is the maximal weighted degree of a monomial in $P$. For the zero polynomial we define the weighted degree as $-\infty$.
The weighted degree of a vector of polynomials is the maximal weighted degree of a polynomial in the vector.

Fix the elements $b_1,\ldots,b_t$ from \cref{claim:unique-decomp-b_i}.
For $s\in[1,d]$, let $M_s\in\Z^{r_s\times r_s}$ be the matrix with
rows and columns indexed by $I_s$ and entry $\Hall_i(b_j)$ in row $i$
and column $j$.

\begin{lemma}
\label{lem:compute-M_s}
Every integer matrix $M_s$ can be computed in polynomial time
and is invertible over $\Q$.
\end{lemma}

\begin{proof}
\cref{claim:unique-decomp-b_i} yields
 $\Hall_i(b_j)=0$ for all $s\in[1,d]$, $i,j\in I_s$ with $i<j$ and $\Hall_i(b_i)=e_i \neq 0$.
Thus $M_s$ is lower triangular with non-zero diagonal entries and therefore invertible over $\Q$.
By Lemmas~\ref{lem:slp-hall-coordinates} and \ref{claim:unique-decomp-b_i}, one can compute all integers $\Hall_i(b_j)$
and hence all matrices $M_s$ in polynomial time.
\end{proof}

\begin{lemma}
\label{lem:triangular-kernel-coordinates1}
One can compute in polynomial time polynomials
$S_i\in\Q[\overline{x}]$ of weighted degree at most $\omega_i$ such that
\begin{equation}
 S_i(\overline{z})=\Hall_i(b_1^{z_1}\cdots b_t^{z_t})
 \label{eq:poly-S_i}
\end{equation}
for every $\overline{z}\in\Z^t$ and $i\in[1,t]$.
\end{lemma}

\begin{proof}
The existence  of the polynomials $S_i$ satisfying \eqref{eq:poly-S_i} follows  
from \cite[Theorem~6.5]{Hall1969}.
We verify the weighted degree bound for the polynomials $S_i$.
For this, we use the Magnus embedding $\mu_{d,r}$.

Fix $i\in[1,t]$ and let
$W_i=\{w\in\{\chi_1,\ldots,\chi_r\}^*:1\leq |w|\leq \omega_i\}$.
By \cite[Corollary~6.13]{Reutenauer1993}, there are integers
$c_{i,w}$ for $w\in W_i$ such that
\begin{equation} \label{eq:hall-from-magnus}
\Hall_i(g)=\sum_{w\in W_i}c_{i,w} \cdot [w]\mu_{d,r}(g)
\end{equation}
for every $g\in\mathcal{N}_{d,r}$.
Reutenauer states this identity for the free group
$\mathcal{F}_r$. Since $\omega_i\leq d$, both the Hall coordinates $\Hall_i(g)$ and
the coefficients $[w]\mu_{d,r}(g)$ depend only
on the image of $g$ in
$\mathcal{F}_r/\gamma_{d+1}(\mathcal{F}_r)=\mathcal{N}_{d,r}$.

For $j\in[1,t]$ let
$U_j=\mu_{d,r}(b_j)-1\in\mathcal{A}_{d,r}$.
Since $b_j\in N_j\subseteq\gamma_{\omega_j}(\mathcal{N}_{d,r})$,
every monomial of $U_j$ has length at least $\omega_j$ \cite[Lemma~6.10(iii)]{Reutenauer1993}.
Consequently, $U_j^q=0$ whenever $q\cdot \omega_j>d$, and the finite
binomial expansion gives for every $z\in\Z$, 
\[\mu_{d,r}(b_j^z)=
\sum_{q=0}^{\lfloor d/\omega_j\rfloor}\binom{z}{q}U_j^q.\]
Let $w\in W_i$. Expanding the product 
\[ \mu_{d,r}\bigg(\prod_{j=1}^t b_j^{x_j}\bigg) = \prod_{j=1}^t \mu_{d,r}(b_j^{x_j}) = \prod_{j=1}^t \sum_{q=0}^{\lfloor d/\omega_j\rfloor}\binom{x_j}{q}U_j^q
\]
shows that
$[w]\mu_{d,r}(b_1^{x_1}\cdots b_t^{x_t})$
is a polynomial $F_w\in\Q[\overline{x}]$. More precisely, it is a sum of
 terms of the form
\begin{equation} \label{eq:weight-bound}
\bigl([w](U_1^{q_1}\cdots U_t^{q_t})\bigr)
\prod_{j=1}^t\binom{x_j}{q_j},
\end{equation}
where $\overline{q}=(q_1,\ldots,q_t)\in\N^t$.
Every monomial of $U_1^{q_1}\cdots U_t^{q_t}$ has length at least
$\sum_{j=1}^t q_j\omega_j$. Hence the coefficient
$[w](U_1^{q_1}\cdots U_t^{q_t})$ can be nonzero only if
$\sum_{j=1}^t q_j\omega_j\leq |w|$.
Since $x_j$ has weight $\omega_j$ and $\binom{x_j}{q_j}$ has ordinary
degree $q_j$, a non-zero term of the form \eqref{eq:weight-bound} has weighted degree $\sum_{j=1}^t q_j\omega_j\leq |w|$.
Thus $F_w$ has weighted degree at most $|w|$.

Applying \eqref{eq:hall-from-magnus} to
$g=b_1^{z_1}\cdots b_t^{z_t}$ yields
\begin{equation} \label{eq:S_i-final}
S_i(\overline{z}) =  \Hall_i(b_1^{z_1}\cdots b_t^{z_t})=
\sum_{w\in W_i}c_{i,w} \cdot F_w(\overline{z})
\end{equation}
for every $\overline{z}\in\Z^t$.
Since polynomials over $\Q$ that agree on $\Z^t$ are identical,
this yields an identity in $\Q[\overline{x}]$.
Every summand on the right-hand side of \eqref{eq:S_i-final}
has weighted degree at most $|w|\leq \omega_i$.
Therefore $S_i$ has weighted degree at most $\omega_i$.

We finally show that every polynomial $S_i$ can be computed in polynomial time. Fix an $i \in [1,t]$ and
consider the set $\mathcal{D}_{d,t} = \{ \overline{\delta}\in\N^t : |\overline{\delta}|_1\leq d \}$. It consists 
of $\binom{t+d}{d} \leq \mathcal{O}(t^d) = \mathcal{O}\big(m^{d^2}\big)$ tuples. 
By \cite[Corollary~2]{ChungYao1977} the set $\mathcal{D}_{d,t}$ admits unique multivariate polynomial interpolation for polynomials in $t$ variables
and of degree at most $d$. To determine $S_i$, we compute in polynomial time $\Hall_i(b_1^{\delta_1}\cdots b_t^{\delta_t}) = S_i(\overline{\delta})$
for every $\overline{\delta} \in \mathcal{D}_{d,t}$ using Lemmas~\ref{lem:slp-hall-coordinates} and \ref{claim:unique-decomp-b_i}. Then $S_i$ is the unique polynomial $P$ in $t$ variables and of degree at most $d$ such that
$S_i(\overline{\delta})=P(\overline{\delta})$ for all  $\overline{\delta} \in \mathcal{D}_{d,t}$. 
Since $S_i$ has at most $\binom{t+d}{d}$ monomials, we can write down for every $\overline{\delta} \in \mathcal{D}_{d,t}$ a linear equation 
in variables $x_{\overline{\alpha}}$ for $\overline{\alpha} \in \mathcal{D}_{d,t}$:
\[  \sum_{\overline{\alpha} \in \mathcal{D}_{d,t}}  \big( \prod_{i=1}^t \delta_i^{\alpha_i} \big) \cdot x_{\overline{\alpha}} = S_i(\overline{\delta}).\]
The matrix corresponding to this system of linear equations must be non-singular. Moreover, every matrix entry has bit length $\mathcal{O}(d \log d)$ 
and every right-hand side $S_i(\overline{\delta})$ also has polynomial bit length (since it has been computed in polynomial time). 
Hence, we can compute the unique solution of the system in polynomial time using Gaussian elimination \cite[Theorem~3.3]{Schrijver1986}, which yields
all coefficients of $S_i$. 
\end{proof}

\begin{lemma}
\label{lem:triangular-kernel-coordinates2}
One can compute in polynomial time column vectors
$\overline{R}_s\in\Q[\overline{x}_{<s}]^{r_s}$ of weighted degree at most $s$
such that
\begin{equation}
 \bigl(S_i(\overline{x})\bigr)^{\mathsf T}_{i\in I_s}
 =M_s\cdot\overline{x}^{\mathsf T}_s+\overline{R}_s(\overline{x}_{<s}).
 \label{eq:triangular-hall-coordinates}
\end{equation}
Moreover, $\overline{R}_s(\overline{0})=\overline{0}$ and
$\overline{R}_s(\Z^{|I_{<s}|})\subseteq\Z^{r_s}$ for every $s\in[1,d]$.
\end{lemma}

\begin{proof}
Fix $s\in[1,d]$ and $i\in I_s$.
Since the polynomial $S_i$ has weighted degree at most $s$ (\cref{lem:triangular-kernel-coordinates1}) and every
variable $x_j$ has weight $\omega_j \ge 1$,  a variable $x_j$
with $\omega_j>s$ cannot occur in $S_i$. If $\omega_j=s$ then the only
possible monomial containing $x_j$ is $x_j$ itself: a higher
power or a product with another variable would result in a weighted
degree greater than $s$. All remaining monomials involve only
variables of weight less than $s$. Thus, we can write $S_i$ as
\[ S_i(\overline{x})=\sum_{j\in I_s}c_{ij}x_j+R_i(\overline{x}_{<s})\]
for coefficients $c_{ij}\in\Q$ ($j\in I_s$) and a polynomial
$R_i\in\Q[\overline{x}_{<s}]$ of weighted degree at most $s$.
Let $\overline{R} = (R_1, \ldots, R_t)^{\mathsf T}$.

We now determine the coefficients $c_{ij}$. 
All Hall coordinates of $b_1^{0}\cdots b_t^{0}=1$ are $0$, so $S_i(\overline{0})=0$ and hence
$R_i(\overline{0})=0$. In particular, $\overline{R}_1= \overline{0}$.
For $j\in I_s$, let $\overline{e}_j\in\Z^t$
be the $j$-th standard basis vector. Then
$c_{ij}=S_i(\overline{e}_j) = \Hall_i(b_j)$.
In particular, $c_{ij}\in\Z$, and
$(c_{ij})_{i,j\in I_s}=M_s$ is the matrix from \cref{lem:compute-M_s}. This proves
\eqref{eq:triangular-hall-coordinates}.

Finally, take an integer vector $\overline{z} = (z_1, \ldots, z_t) \in \Z^t$ with $\overline{z}_s = \overline{0}$.
Then \eqref{eq:triangular-hall-coordinates} yields
$\overline{R}_s(\overline{z}_{<s}) = \bigl(\Hall_i(\prod_{j\in I_{<s}}b_j^{z_j})\bigr)^{\mathsf T}_{i\in I_s} \in \Z^{|I_s|}$.
Hence, we have $\overline{R}_s(\Z^{|I_{<s}|})\subseteq\Z^{r_s}$.

Once the polynomials $S_i$ have been computed (\cref{lem:triangular-kernel-coordinates1}), the equations \eqref{eq:triangular-hall-coordinates} allow us to determine all vectors $\overline{R}_s$ in polynomial time.
\end{proof}
Introduce now new commuting variables $\overline{y}=(y_1,\ldots,y_m)$ and let
\[\overline{v}=(y_1,\ldots,y_m,-1,0,\ldots,0)^{\mathsf{T}}\in\Z[\overline{y}]^t.\]
For an integer vector $\overline{y}\in\Z^m$, this is the vector of 
 Hall coordinates of $a_1^{y_1}\cdots a_m^{y_m}a_r^{-1}$. 
Note that $\overline{v}_1 = (y_1,\ldots,y_m,-1)^{\mathsf{T}}$ and $\overline{v}_s = \overline{0}^{\mathsf{T}}$ for $s \in [2,d]$. 

By induction on $s \in [1,d]$ we define the column vector of polynomials $\overline{Q} = (Q_1,\ldots,Q_t)^{\mathsf{T}} \in\Q[\overline{y}]^t$ as follows:
\begin{equation}
 \overline{Q}_s=M_s^{-1}\cdot \bigl(\overline{v}_s-\overline{R}_s(\overline{Q}_{<s})\bigr).
 \label{eq:kernel-coordinate-recursion}
\end{equation}
\begin{lemma}
\label{eq:kernel-integrality}
The polynomials $Q_i$ can be computed in polynomial time and satisfy
$\deg Q_i\leq \omega_i$ for $i\in[1,t]$. For every
$\overline{\epsilon}=(\epsilon_1,\ldots,\epsilon_m)\in\{0,1\}^m$,
\begin{equation}
 g_1^{\epsilon_1}\cdots g_m^{\epsilon_m}=h \text{ in } H
 \quad\Longleftrightarrow\quad
 \forall i\in[1,t]:Q_i(\overline{\epsilon})\in\Z.
 \label{eq:polynomial-integrality-test}
\end{equation}
\end{lemma}

\begin{proof}
The weighted degree bound $s$ on the polynomials in $\overline{R}_s$
implies inductively that all polynomials in $\overline{Q}_s$ have
ordinary degree at most $s \leq d$. Indeed, a monomial
$\prod_{j\in I_{<s}}x_j^{\alpha_j}$ of weighted degree at most $s$ becomes
$\prod_{j\in I_{<s}}Q_j^{\alpha_j}$ of ordinary degree at most
$\sum_{j\in I_{<s}}\alpha_j\cdot\omega_j\leq s$.
The entries of $\overline{v}_1$ have degree at most one,
$\overline{v}_s=\overline{0}$ for $s\geq2$ and multiplication with the matrix $M^{-1}_s$ does not change the degree.

By Lemmas~\ref{lem:compute-M_s} and \ref{lem:triangular-kernel-coordinates2}, all $M_s$ and $\overline{R}_s$
are computable in polynomial time. All inverses $M_s^{-1}$ can also be
computed in polynomial time \cite[Corollary~3.3a(iii)]{Schrijver1986}.
The $d=\mathcal{O}(1)$ rounds of substitution of $\overline{Q}_{<s}$
into $\overline{R}_s$ and multiplication by $M_s^{-1}$ in
\eqref{eq:kernel-coordinate-recursion} compute every $\overline{Q}_s$
in polynomial time and preserve polynomial bit length of all coefficients.
The degree bound ensures that each intermediate polynomial has only
polynomially many monomials.

It remains to prove \eqref{eq:polynomial-integrality-test}.
Since $\overline{v}_s = M_s\cdot\overline{Q}_s + \overline{R}_s(\overline{Q}_{<s}) = (S_i(\overline{Q}))^{\mathsf T}_{i \in I_s}$ for all $s \in [1,d]$,
it follows that for every $\overline{\epsilon}\in\{0,1\}^m$, the vector
$\overline{Q}(\overline{\epsilon}) \in \Q^t$ is the unique solution of the system
$\bigl(S_i(\overline{x})\bigr)^{\mathsf T}_{i\in[1,t]}=\overline{v}(\overline{\epsilon})$
(uniqueness follows from the fact that every $M_s$ is invertible).
If $\overline{Q}(\overline{\epsilon}) \in \Z^t$ then we obtain
\[ b_1^{Q_1(\overline{\epsilon})}\cdots b_t^{Q_t(\overline{\epsilon})}=a_1^{\epsilon_1}\cdots a_m^{\epsilon_m}a_r^{-1}\]
in $\mathcal{N}_{d,r}$.
Since the left-hand side belongs to the kernel $K$, the same holds for the right-hand side, which yields $g_1^{\epsilon_1}\cdots g_m^{\epsilon_m}=h$ in $H$.

Conversely, if $g_1^{\epsilon_1}\cdots g_m^{\epsilon_m}=h$ then 
$a_1^{\epsilon_1}\cdots a_m^{\epsilon_m}a_{r}^{-1} \in K$. Hence,
\cref{claim:unique-decomp-b_i} gives integers $z_1, \ldots, z_t \in \Z$ with 
\[ b_1^{z_1}\cdots b_t^{z_t} = a_1^{\epsilon_1}\cdots a_m^{\epsilon_m}a_{r}^{-1}. 
\]
This yields an integer 
solution of the system $\bigl(S_i(\overline{x})\bigr)^{\mathsf T}_{i\in[1,t]}=\overline{v}(\overline{\epsilon})$. Since this system has the unique solution $\overline{Q}(\overline{\epsilon})$, we get
$\overline{Q}(\overline{\epsilon}) \in \Z^t$.
\end{proof}
Let $\rho=\operatorname{lcm}\{1,\ldots,n\}\in\N$.

\begin{lemma}
\label{denom-bound-M_s}
For every $s\in[1,d]$, $\rho\cdot M_s^{-1}\in\Z^{r_s\times r_s}$.
\end{lemma}

\begin{proof}
Let $s \in [1,d]$. Recall that  $g \in \mathcal{N}_{d,r}$ belongs to $\gamma_s(\mathcal{N}_{d,r})$ if 
and only if $\Hall_i(g) = 0$ for all $i \in I_{<s}$.
Moreover, the cosets $a_i \gamma_{s+1}(\mathcal{N}_{d,r})$ with $i \in I_s$
form a free abelian basis of the free abelian group 
$\gamma_s(\mathcal{N}_{d,r})/\gamma_{s+1}(\mathcal{N}_{d,r}) \cong \mathbb{Z}^{r_s}$; see 
\cite[Corollary~6.16]{Reutenauer1993}. An isomorphism from $\gamma_s(\mathcal{N}_{d,r})/\gamma_{s+1}(\mathcal{N}_{d,r})$ to $\mathbb{Z}^{r_s}$
is obtained by mapping $g\gamma_{s+1}(\mathcal{N}_{d,r})$, for $g\in\gamma_s(\mathcal{N}_{d,r})$, to
$(\Hall_i(g))^{\mathsf T}_{i\in I_s}$ \cite[Corollary~6.16]{Reutenauer1993}. 

We claim that under the corresponding projection from $\gamma_s(\mathcal{N}_{d,r})$ to $\mathbb{Z}^{r_s}$, the image
of $K\cap\gamma_s(\mathcal{N}_{d,r})$ is the lattice $M_s\Z^{r_s}$.
Let $g\in K\cap\gamma_s(\mathcal{N}_{d,r})$ have the factorization
$g=b_1^{\alpha_1}\cdots b_t^{\alpha_t}$ from \cref{claim:unique-decomp-b_i},
where $\overline{\alpha} = (\alpha_1,\ldots,\alpha_t)^{\mathsf T}\in\Z^t$.
Then $\Hall_i(g) = 0$ for $i \in I_{<s}$.
Using \eqref{eq:triangular-hall-coordinates}, the invertibility
of the matrices $M_1,\ldots,M_{s-1}$, and
$\overline{R}_\ell(\overline{0})=\overline{0}$ for every $\ell\in[1,s-1]$,
we obtain successively $\overline{\alpha}_{<s} = \overline{0}$.
Hence $(\Hall_i(g))^{\mathsf T}_{i\in I_s}=M_s \overline{\alpha}_{s}$ by \eqref{eq:triangular-hall-coordinates}.
Conversely, for every integer vector $\overline{\alpha}_s = (\alpha_i)^{\mathsf T}_{i\in I_s}\in\Z^{r_s}$
one has $g = \prod_{i\in I_s}b_i^{\alpha_i} \in K\cap\gamma_s(\mathcal{N}_{d,r})$ with
$(\Hall_i(g))^{\mathsf T}_{i\in I_s} = M_s\overline{\alpha}_s$.

 Since $H \leq \Sym(n)$,
every element of $H$ has order dividing $\rho$. Hence,
for every $i\in I_s$, the element $a_i^\rho$ belongs to $K  \cap\gamma_s(\mathcal{N}_{d,r})$.
The vector $(\Hall_j(a_i^\rho))^{\mathsf T}_{j\in I_s}$ has entry $\rho$ at $i$
and zero elsewhere. 
Each of these vectors belongs to $M_s \cdot\Z^{r_s}$ by the previous paragraph. 
It  follows
that $\rho\cdot\Z^{r_s}\subseteq M_s \cdot\Z^{r_s}$, i.e., 
$\rho \cdot M_s^{-1}\in\Z^{r_s\times r_s}$.
\end{proof}
Define $\alpha_1=1$ and $\alpha_s=1+d\alpha_{s-1}\in\N$ for
$s\in[2,d]$ and let $\beta=d!\rho\in\N$ and $q = \beta^{\alpha_d}$. Note that $\alpha_d$ only
depends on $d$.

\begin{lemma}
\label{claim:denominators1}
For every $s\in[1,d]$, $\beta^{\alpha_s}\cdot\overline{Q}_s\in\Z[\overline{y}]^{r_s}$.
\end{lemma}

\begin{proof}
We first recall a fact about integer-valued polynomials.
For $\ell \geq1$, consider a polynomial $P\in\Q[x_1,\ldots,x_\ell]$ of degree at most $d$
 such that $P(\Z^\ell) \subseteq \Z$.
Then, $P$ can be written as 
$P = \sum_{\overline{z} \in C} c_{\overline{z}}\prod_{i=1}^\ell\binom{x_i}{z_i}$
where $C = \{ \overline{z} = (z_1, \ldots, z_\ell) \in \N^\ell : \sum_{j=1}^\ell z_j\leq d \}$ and
$c_{\overline{z}}\in\Z$ for all $\overline{z} \in C$; see e.g.~\cite[Section~2.4]{ChabertCahen2006} and \cite[pp.~121--123]{Ostrowski1919}.
Since $\prod_{i=1}^\ell z_i!$ divides $(\sum_{i=1}^\ell z_i)!$ (the quotient is a multinomial coefficient) 
it follows that $d! \cdot P \in \Z[x_1,\ldots,x_\ell]$.

Every polynomial $R_i$ with $i\in I_s$ and $s\in[1,d]$ has
 ordinary degree at most $s\leq d$,
since every variable in $\overline{x}_{<s}$ has 
weight at least one and $R_i$ has weighted degree at most $s$ by \cref{lem:triangular-kernel-coordinates2}.
Applying the observation from the previous paragraph shows 
that $d!\overline{R}_s\in\Z[\overline{x}_{<s}]^{r_s}$.
For $s=1$, the conclusion follows directly from $\overline{R}_1=\overline{0}$.

We show $\beta^{\alpha_s}\overline{Q}_s\in\Z[\overline{y}]^{r_s}$
by induction on $s$.
We have $\overline{Q}_1=M_1^{-1}(y_1,\ldots,y_m,-1)^{\mathsf T}$ by \eqref{eq:kernel-coordinate-recursion} and $\overline{R}_1=\overline{0}$.
By \cref{denom-bound-M_s}, $\rho\overline{Q}_1\in\Z[\overline{y}]^r$
and hence $\beta^{\alpha_1}\overline{Q}_1\in\Z[\overline{y}]^r$.
For $s\in[2,d]$, we have
$\overline{Q}_s=-M_s^{-1} \cdot \overline{R}_s(\overline{Q}_{<s})$ by \eqref{eq:kernel-coordinate-recursion}.
Since the sequence $(\alpha_s)_{s\in[1,d]}$ is increasing, the induction
hypothesis shows $\beta^{\alpha_{s-1}} \overline{Q}_{<s} \in\Z[\overline{y}]^{|I_{<s}|}$.
Each polynomial in $\overline{R}_s$ has ordinary degree at most $d$, and
$d!\overline{R}_s$ has integer coefficients. Applying these facts to
each monomial of $\overline{R}_s(\overline{Q}_{<s})$ gives
$d!\beta^{d\alpha_{s-1}}\overline{R}_s(\overline{Q}_{<s})\in\Z[\overline{y}]^{r_s}$.
Finally, with \cref{denom-bound-M_s} we get
$\beta^{\alpha_s}\overline{Q}_s
=-\rho M_s^{-1}\cdot d!\beta^{d\alpha_{s-1}}
\overline{R}_s(\overline{Q}_{<s})\in\Z[\overline{y}]^{r_s}$.
\end{proof}

For a prime $p$, write $v_p(q)\ge 0$ for the largest $v \in \N$ such that $p^v$ divides $q$.
\begin{lemma}
\label{claim:denominators2}
The integer $q$ and its prime factorization can be computed
in time $\poly(n)$. For every prime divisor $p$ of $q$, 
$p^{v_p(q)} \leq \poly(n)$ and the number of prime divisors of $q$ is also bounded by $\poly(n)$.
\end{lemma}

\begin{proof} 
To factorize $q = (d! \rho)^{\alpha_d}$ it suffices to factorize $\rho$ and $d!$.
The latter is a fixed constant. To factorize $\rho = \operatorname{lcm}\{1,\ldots,n\}$ it suffices
to factorize all numbers  in $[2,n]$, which can be done in time $\poly(n)$.
If $p$ divides $q$ then $p^{v_p(q)}
 =p^{\alpha_d(v_p(d!)+v_p(\rho))}
 \leq (d!)^{\alpha_d}n^{\alpha_d} = \poly(n)$, 
since $p^{v_p(\rho)}\leq n$, also when $p>n$.
Since $\log\rho\leq\log(n!)=\poly(n)$ and $\alpha_d$ only depends on
$d$, $q$ has polynomial bit length. This also shows that $q$ can be computed in polynomial 
time using arithmetic on numbers of bit length $\poly(n)$. The number of distinct
prime factors is bounded by $\log_2q = \poly(n)$.
\end{proof}
We can now conclude the proof of \cref{lem:nilpotent-polynomials}.

\begin{proof}[Proof of \cref{lem:nilpotent-polynomials}]
Construct the kernel basis from \cref{claim:unique-decomp-b_i} and the
polynomials $Q_1,\ldots,Q_t$ from \cref{eq:kernel-integrality} and set
$P_i=q\cdot Q_i$ for $i\in[1,t]$. By \cref{claim:denominators1}, $P_i\in\Z[\overline{y}]$.
Moreover, $\deg P_i\leq \omega_i\leq d$ and $t=\mathcal{O}(m^d)$.
For every $\overline{\epsilon}\in\{0,1\}^m$, the condition
$P_i(\overline{\epsilon})\equiv0\pmod q$ is equivalent to
$Q_i(\overline{\epsilon})\in\Z$.
Thus \cref{eq:kernel-integrality} proves \eqref{eq:nilpotent-polynomials}.
All constructions take polynomial time, and \cref{claim:denominators2}
provides the required factorization and prime-power bounds for $q$.
\end{proof}

\subsection{Proof of the main result on nilpotent permutation groups}

Using \cref{lem:nilpotent-polynomials} we can prove our main result for nilpotent permutation groups:
\begin{theorem}
\label{thm:nilpotent}
For every fixed $d \geq 1$, there is an fpt-reduction
from $\pSSS[\NilPG{d}]$ to $\pSSS[\SolPG{1}]$.
Consequently, $\pSSS[\NilPG{d}]$ belongs to $\W[2]$.
\end{theorem}

\begin{proof}
The case $d=1$ is trivial since $\pSSS[\NilPG{1}]$ and $\pSSS[\SolPG{1}]$ are the same problem. Assume now $d \ge 2$.
Consider an instance consisting of 
$g_1,\ldots,g_m,h\in \Sym(n)$ with $G := \langle g_1,\ldots,g_m\rangle$ $d$-step nilpotent and parameter $k\in [1,m]$. W.l.o.g.~we can assume that $h \in G$.
Compute $t,q\ge 1$ and
$P_1,\ldots,P_t\in\Z[y_1,\ldots,y_m]$ as in
\cref{lem:nilpotent-polynomials}.

Let $A=(\Z/q\Z)^t$, written additively. Replace in every $P_i$ all non-zero exponents of variables $y_j$ by $1$.
This neither changes values on $\{0,1\}^m$ nor increases the degree. After reducing coefficients modulo $q$, write
\begin{equation}
 (P_1(\overline{y}),\ldots,P_t(\overline{y}))
 =\sum_{\substack{S\subseteq[1,m]\\|S|\leq d}}
 v_S\prod_{i\in S}y_i
 \quad\text{with}\quad v_S\in A.
 \label{eq:multilinear-coefficients}
\end{equation}
Let $u=-v_{\emptyset}\in A$.
For a set $I\subseteq[1,m]$ of size $k$, the equivalence
\eqref{eq:nilpotent-polynomials} and
\eqref{eq:multilinear-coefficients} show that $\prod_{i \in I} g_i = h$,
where the product is taken in increasing order of $i$, if and only if
\begin{equation}
 \sum_{\substack{\emptyset\neq S\subseteq I\\|S|\leq d}}v_S=u.
 \label{eq:subset-coefficient-sum}
\end{equation}
Let $\mathcal T=\{P\subseteq[1,k]:1\leq|P|\leq d\}$ and
$\ell=|\mathcal T|=\sum_{s=1}^{\min\{d,k\}}\binom{k}{s} \le k^{d+1}$.
The sets in $\mathcal T$ will be called types.
Let $\mathcal E=\{(P,p):P\in\mathcal T,\ |P|\geq2,\ p\in P\}$.
We construct an instance of $\pSSS[\SolPG{1}]$ over the abelian group
\begin{equation}
 A' =A\times(\Z/(\ell+1)\Z)^{\mathcal T}
          \times(\Z/(m+1)\Z)^{\mathcal E}
 \label{eq:nilpotent-abelian-output}
\end{equation}
with target $(u,\overline{1},\overline{0})\in A'$ and parameter $\ell$.
Here $\overline{1}$ has entry $1$ at each coordinate from $\mathcal T$,
and $\overline{0}$ is zero at each coordinate from $\mathcal E$.

For each $P\in\mathcal T$ and each injective function $\iota : P \to [1,m]$ produce one output element $a[P,\iota] \in A'$.
It is defined as $a[P,\iota] = (v_{\iota(P)}, T, E)$,
where $v_{\iota(P)} \in A$ is from \eqref{eq:multilinear-coefficients} and the mappings $T : \mathcal T \to \Z/(\ell+1)\Z$ and
$E : \mathcal E \to \Z/(m+1)\Z$ are defined as follows:
\begin{itemize}
\item $T(P) = 1$ and $T(P') = 0$ for all $P' \in \mathcal{T} \setminus \{P\}$.
\item If $P=\{p\}$ is a singleton then $E(P',p) = \iota(p)$ for every $P' \in \mathcal{T}$ with $|P'| \ge 2$ and $p \in P'$. 
If $|P|\geq2$ then $E(P,p) = -\iota(p)$ for every $p\in P$. All other values of $E$ are $0$.
\end{itemize}
The order among the elements $a[P,\iota]$ in the output list is arbitrary since $A'$
is abelian. 

Suppose that a subset $S$ of $\ell$ output elements $a[P,\iota]$
has sum $(u,\overline{1},\overline{0})\in A'$. 
The $(\Z/(\ell+1)\Z)^{\mathcal T}$-component enforces that for every type $P$, the set $S$
contains exactly one element of the form $a[P,\iota_P]$, i.e., 
 $S = \{a[P,\iota_P]:P\in\mathcal T\}$.
For a type $P\in\mathcal T$ with $|P|\geq2$ and $p \in P$, the coordinate $(P,p)$ of $(\Z/(m+1)\Z)^{\mathcal E}$ forces
$\iota_{\{p\}}(p) - \iota_P(p) \equiv0\pmod{m+1}$. 
Since both indices belong to $[1,m]$, this implies
\begin{equation} \label{eq:iota_P}
\iota_{\{p\}}(p) = \iota_P(p).
\end{equation}
Define $\iota:[1,k]\to[1,m]$ by $\iota(p)=\iota_{\{p\}}(p)$ for $p\in[1,k]$.
Equality \eqref{eq:iota_P} implies $\iota_P=\iota\rest_P$ for every $P\in\mathcal T$.
If $k\geq2$ then $\{p,p'\}\in\mathcal T$ for all distinct $p,p'\in[1,k]$, since $d\geq2$.
Injectivity of $\iota_{\{p,p'\}}$ therefore gives $\iota(p)\neq\iota(p')$.
For $k=1$, the map $\iota$ is injective as well. Hence $I=\iota([1,k])\subseteq[1,m]$ has size $k$.
As $P$ ranges over $\mathcal T$, the sets $\iota(P)$ run repetition-free through all nonempty subsets of $I$ of size at most $d$.
Thus the $A$-components $v_{\iota(P)}$ of the elements $a[P,\iota_P]$ sum to the left-hand side of
\eqref{eq:subset-coefficient-sum}. On the other hand, the sum of all $a[P,\iota_P]$ ($P\in\mathcal T$) is
the target $(u,\overline{1},\overline{0})$. Hence, \eqref{eq:subset-coefficient-sum} holds,
and therefore $\prod_{i \in I} g_i = h$.

Conversely, let $I\subseteq[1,m]$ be such that $|I|=k$ and $\prod_{i \in I} g_i = h$ holds.
This implies that \eqref{eq:subset-coefficient-sum} holds.
Choose a bijection $\iota:[1,k]\to I$. Let $S = \{ a[P,\iota\rest_P] : P \in \mathcal{T}\}$.
We have $|S|=\ell$. Moreover, the $(\Z/(\ell+1)\Z)^{\mathcal T}$-components in $S$
sum to $\overline{1}$.
At each coordinate $(P,p)\in\mathcal E$ of $(\Z/(m+1)\Z)^{\mathcal E}$ (note that $|P| \geq 2$ and $p \in P$),
the element $a[\{p\},\iota\rest_{\{p\}}]$ contributes $\iota(p)$,
the element $a[P,\iota\rest_P]$  contributes $-\iota(p)$, and all other elements in $S$ contribute zero.
Hence the $(\Z/(m+1)\Z)^{\mathcal E}$-components  in $S$ sum to $\overline{0}$. Finally, the
$A$-components in $S$ sum to $u$ by \eqref{eq:subset-coefficient-sum}, since the mapping
$P\mapsto\iota(P)$ ($P \in \mathcal{T}$) enumerates all nonempty subsets of $I$ of size at most $d$ without repetitions.
This proves the correctness of our reduction.

It remains to realize the group $A'$ as a permutation group on a small domain.
Let $q = p_1^{v_1} \cdots p_\lambda^{v_\lambda}$ be the prime factorization of $q$. The Chinese remainder theorem identifies $\Z/q\Z$ with
$\prod_{i=1}^\lambda \Z/p_i^{v_i}\Z$.

A cyclic group $\Z/a\Z$ can be represented as a cyclic subgroup of $\Sym(a)$.
Moreover, the direct product $G \times H$ of permutation groups $G \leq \Sym(a)$ and $H \leq \Sym(b)$
can be represented as a subgroup of $\Sym(a+b)$. 
Hence, our abelian group $A'$ can be represented as a permutation group on a domain of size
\begin{equation}
n' :=  t\sum_{i=1}^\lambda p_i^{v_i}
   +\ell(\ell+1)+|\mathcal E|(m+1).
 \label{eq:nilpotent-output-degree}
\end{equation}
The output elements $a[P,\iota] \in A'$ as well as the target $(u,\overline{1},\overline{0})$
can be realized as elements of $\Sym(n')$.
By \cref{lem:nilpotent-polynomials}, $t = \poly(m)$,
$\lambda = \poly(n)$, and
$p_i^{v_i} = \poly(n)$ for all $i \in [1,\lambda]$. Recall that $\ell = |\mathcal{T}| \leq k^{d+1} = \poly(k)$ and hence $|\mathcal E|\leq d\ell = \poly(k)$. With \eqref{eq:nilpotent-output-degree} we obtain $n' = \poly(k,m,n)$.

For each type $P\in\mathcal T$, there are at most $m^{|P|}\leq m^d$
injective functions $\iota:P\to[1,m]$, so the number of output elements is at most $\ell m^d = \poly(m,k)$. Therefore all output permutations can be
produced in time $\poly(m,n,k)$. The output parameter $\ell$ depends only on $d$ (fixed) and $k$ (the parameter of the input instance).
This shows that our reduction is an fpt-reduction. Finally, \cref{thm:solvable} with
$d=1$ yields membership in $\W[2]$.
\end{proof}
\cref{pre:lem:replication,pre:lem:padding} show that for every $d \geq 1$ the problems 
$\pSSS_\le[\NilPG{d}]$, $\pF_\le[\NilPG{d}]$ and $\pF[\NilPG{d}]$ belong to
$\W[2]$ as well.

\cref{thm:nilpotent} states that for every fixed $d$, $\pSSS[\NilPG{d}]$ is in $\W[2]$.
Whether $\pSSS$ for the class of all nilpotent permutation groups is in $\W[2]$ remains open.

\section{W[2]-hardness for nilpotent permutation groups}
\label{sec:nilpotent-hardness}

Let $\pF_{\leq}[\mathsf{NilPG}]$ be the restriction of the problem $\pF_{\leq}[\Sym(*)]$, where
the input permutations $g_1, \ldots, g_m$ generate a nilpotent permutation group; the problems 
$\pF[\mathsf{NilPG}]$, $\pSSS_{\leq}[\mathsf{NilPG}]$ and $\pSSS[\mathsf{NilPG}]$ are defined analogously.
We prove that $\pF_{\leq}[\mathsf{NilPG}]$ is $\W[2]$-hard by a reduction from parameterized \textsc{Set Cover}. 

Fix an integer $k\ge 1$ and let $q=2^{k+1}$. We consider $\Sym([0,q-1])$, where $[0,q-1]$ is identified with 
$\Z/q\Z$, and define two permutations 
$r,s\in\Sym([0,q-1])$ by $x^r=x+1$ and $x^s=-x$ for
$x\in\Z/q\Z$. Let $D_q=\langle r,s\rangle\leq\Sym([0,q-1])$ be the dihedral group for a regular $q$-gon.
Then $r^q=s^2=1$ and $s^{-1}rs=r^{-1}$. The $2q$ maps
$r^i s^\epsilon$ with $i\in[0,q-1]$ and
$\epsilon\in\{0,1\}$ are distinct and comprise $D_q$.
Thus $D_q$ is a $2$-group of order $2q = 2^{k+2}$, represented as a permutation group on
$q = 2^{k+1}$ points. 

The following lemma specializes the commutator construction for conjunctions
used by Wei{\ss} in the proof of \cite[Proposition~8]{Weiss2020}.
The technique goes back to Barrington \cite{Barrington1989} and is also developed
in \cite{BartholdiFigeliusLohreyWeiss2022}.
Here the dihedral group depends on $k$, and we give the proof to obtain the
explicit length and output element needed in the proof of
\cref{thm:nilpotent-hardness}.

\begin{lemma}
\label{lem:dihedral-test}
Given $k\ge 1$, one can construct, in time depending
only on $k$, a sequence $(k_i,g_{i,0},g_{i,1})_{i=1}^{\ell}$
with $\ell=3\cdot2^k-2\in\N$, $k_i\in[1,k]$, and
$g_{i,0},g_{i,1}\in D_q$ for $i\in[1,\ell]$, such that
for every $\overline{x}=(x_1,\ldots,x_k)\in\{0,1\}^k$ the following holds in $D_q$:
\begin{equation}
 \prod_{i=1}^{\ell}g_{i,x_{k_i}}
 =\begin{cases}
   r^{2^k} &\text{if } x_1=\cdots=x_k=0,\\
   1 &\text{otherwise}.
  \end{cases}
 \label{eq:dihedral-test}
\end{equation}
\end{lemma}

\begin{proof}
For $z\in\Z$ and $\epsilon\in\{0,1\}$ one easily shows $[r^z,s^\epsilon]=r^{-2z\epsilon}$ in $D_q$.
Define maps $C_j:\{0,1\}^k\to D_q$ for $j\in[0,k]$ recursively
by $C_0(\overline{x})=r$ and
$C_j(\overline{x})=[C_{j-1}(\overline{x}),s^{1-x_j}]$ for $j\in[1,k]$.
Induction gives
\begin{equation}
 C_j(\overline{x})=r^{(-2)^j\prod_{p=1}^j(1-x_p)}
 \qquad(j\in[0,k],\ \overline{x}\in\{0,1\}^k).
 \label{eq:dihedral-commutators}
\end{equation}
Since $r$ has order $2^{k+1}$, its powers $r^{(-2)^k}$ and
$r^{2^k}$ coincide and are not the identity. Hence
\eqref{eq:dihedral-commutators} with $j=k$ gives the two values
required in \eqref{eq:dihedral-test}.

We then expand the nested commutator  $C_j$ into a word $w_j$, where every symbol is
$r$, $r^{-1}$, $s^{1-x_i}$, or $s^{x_i-1}$ for some $i \in [1,k]$. 
The length $\ell_j = |w_j|$ satisfies $\ell_0=1$ and
$\ell_j=2\ell_{j-1}+2$ for $j\in[1,k]$. Thus
$\ell := \ell_k=3\cdot2^k-2$. To obtain the sequence $(k_i,g_{i,0},g_{i,1})_{i=1}^{\ell}$ we replace
the $i$-th symbol $a_i$ in $w_k$ as follows: If $a_i = r^{\pm 1}$, then $a_i$ is replaced by $(1, a_i, a_i)$.
If $a_i = s^{1-x_j}$ then $a_i$ is replaced by $(j, s, 1)$ and if $a_i = s^{x_j-1}$ then $a_i$ is replaced by $(j, s^{-1}, 1)$.
\end{proof}

\begin{theorem}
\label{thm:nilpotent-hardness}
The problem $\pF_{\leq}[\mathsf{NilPG}]$ is $\W[2]$-hard,
already when the input groups are permutation $2$-groups.
\end{theorem}

\begin{proof}
We prove the theorem by a reduction from parameterized \textsc{Set Cover}.
An instance of parameterized \textsc{Set Cover} consists of a
finite universe $U \neq \emptyset$, sets $S_1,\ldots,S_m\subseteq U$ with $m \ge 1$ and a parameter $k\in [1,m]$. The question is whether 
there is a subset $I \subseteq [1,m]$ with $|I| \leq k$ and $U = \bigcup_{i \in I} S_i$. This problem is $\W[2]$-hard
\cite{DowneyFellows1995}.

Fix the parameterized \textsc{Set Cover} instance $U,S_1,\ldots,S_m,k$. 
W.l.o.g.~we can assume that $m \geq 2$.
Clearly, the instance $U,S_1,\ldots,S_m,k$ is positive if and only if there are
indices $i_1,\ldots,i_k\in[1,m]$, with repetitions allowed, such that
\begin{equation}
\bigcup_{j \in [1,k]} S_{i_j} = U.
 \label{eq:hardness-cover}
\end{equation}
We first construct a permutation $2$-group $G$.
Use \cref{lem:dihedral-test} to construct the sequence
$(k_i,g_{i,0},g_{i,1})_{i=1}^{\ell}$, where
$\ell=3\cdot2^k-2$. Put $k'=k+\ell$ and
$\mu=\lceil\log_2m\rceil \geq 1$. Let $C_2=\Z/2\Z$,
and let $\operatorname{code}:[1,m]\to C_2^\mu$ map $i\in[1,m]$
to the binary expansion of $i-1$ with exactly $\mu$ digits.
This map is injective.

Let $a=(1, 2)$ and $b=(1, 3)(2, 4)$ be elements of $\Sym(4)$,
and put $D=\langle a,b\rangle\leq\Sym(4)$.
Since $a$ and $b$ are involutions and $ab = (1,4,2,3)$,
the group $D$ is dihedral of order $8$. In particular,
$ab\neq ba$.
Let $H=D^{k'-1}$ and $\gamma=(ab)_{j\in[1,k'-1]}\in H$.
For $i\in[1,k']$, define $c_i = (c_{i,1}, \ldots, c_{i,k'-1}) \in H$ as follows, where $j \in [1,k'-1]$.
\begin{equation}
 c_{i,j}=\begin{cases}
 a & \text{if } i=j,\\
 b & \text{if } i=j+1,\\
 1 &\text{otherwise}.
 \end{cases}
 \label{eq:hardness-order-coordinates}
\end{equation}
Note that $c_1 c_2\cdots c_{k'}=\gamma$ in $H$.
Consider then the group
\begin{equation}
 G=D_q^U\times C_2^{k'}\times(C_2^\mu)^{\ell}\times H.
 \label{eq:hardness-group}
\end{equation}
Use the natural action of each $D_q$-copy on a separate set
of $2^{k+1}$ points, let each $C_2$-copy act by a
transposition on a separate pair of points, and let each $D$-copy of $H$
act on a separate set of four points. This gives an embedding of $G$ into $\Sym(n)$ for 
\begin{equation}
 n=2^{k+1}|U|+2k'+2\mu\ell+4(k'-1).
 \label{eq:hardness-degree}
\end{equation}
The group in \eqref{eq:hardness-group} is a finite $2$-group and
therefore nilpotent. Its generators can be given explicitly: use
$r$ and $s$ in each $D_q$-copy, the individual transpositions in
the $C_2$-copies, and $a,b$ in each $D$-copy of $H$.

The permutation group $G \le \Sym(n)$ will be the ambient group of our $\pF_{\leq}[\mathsf{NilPG}]$-instance.
We now define the group elements from $G$ that are used in the $\pF_{\leq}[\mathsf{NilPG}]$-instance.
Let $e_1,\ldots,e_{k'}\in C_2^{k'}$ be the standard basis of $C_2^{k'}$.
For $u\in U$ and $i\in[1,m]$, let $\chi_{u,i}\in\{0,1\}$ be $1$ if $u\in S_i$ and $0$ otherwise.
We then construct the following elements of $G$:
\begin{itemize}
\item For $p\in[1,k]$ and $i\in[1,m]$, construct the group element
\begin{equation} \label{eq:vpi}
v_{p,i}=((1)_{u\in U},e_p,\xi_{p,i},c_p)\in G, \text{ where } \xi_{p,i} = (\xi_{p,i,1}, \ldots, \xi_{p,i,\ell}) \in(C_2^\mu)^{\ell}
\end{equation}
with $\xi_{p,i,j} = \operatorname{code}(i)$ whenever $k_j = p$ and 
$\xi_{p,i,j} = \overline{0}$ otherwise.
\item For $t\in[1,\ell]$ and $j\in[1,m]$, construct the group element
\begin{equation} \label{eq:wti}
 w_{t,j}=((g_{t,\chi_{u,j}})_{u\in U},e_{k+t},\eta_{t,j},c_{k+t})\in G, \text{ where } \eta_{t,j} = (\eta_{t,j,1}, \ldots, \eta_{t,j,\ell}) \in(C_2^\mu)^{\ell} 
 \end{equation}
with $\eta_{t,j,t} = \operatorname{code}(j)$ and $\eta_{t,j,i} = \overline{0}$ if $i\neq t$.
\end{itemize}
The target element is 
\[ h=((1)_{u\in U},\overline{1},\overline{0},\gamma)\in G,\] where
$\overline{1}\in C_2^{k'}$ has every entry equal to $1$ and
$\overline{0}\in(C_2^\mu)^{\ell}$ is zero. The output parameter
is $k'$. This concludes the construction of the $\pF_{\leq}[\mathsf{NilPG}]$-instance.
Since $G$ is a permutation $2$-group, the same holds also for the subgroup generated by
all constructed elements $v_{p,i}$ and $w_{t,j}$.

Let us show the correctness of our reduction.
Suppose first that a product of at most $k'$ elements from \eqref{eq:vpi} and \eqref{eq:wti} 
equals $h$, where elements may be repeated and occur in any order.
Each factor contributes $1$ in exactly
one coordinate from $C_2^{k'}$ in \eqref{eq:hardness-group}. The target has a $1$ in all $k'$ 
coordinates. Thus, for every $p \in [1,k]$ there must be an odd number $\alpha_p \geq 1$ of elements of the form
$v_{p,i}$ in the product. Similarly, for every $t \in [1,\ell]$ there must be an odd number $\beta_t \geq 1$ of elements of the form
$w_{t,i}$ in the product. The bound $k' = k+\ell$ on the number of elements in the product
therefore enforces $\alpha_p = \beta_t = 1$ for all $p \in [1,k]$, $t \in [1,\ell]$. Hence, the product consists of $k$ elements 
 $v_{p,i_p}$ for all $p\in[1,k]$ and $\ell$ elements $w_{t,j_t}$ for $t\in[1,\ell]$, where $i_p,j_t\in[1,m]$.

We next show that the product must be 
\begin{equation} \label{eq:order-of-elements}
v_{1,i_1} v_{2,i_2} \cdots v_{k,i_k} w_{1,j_1} w_{2,j_2} \cdots w_{\ell, j_\ell}.
\end{equation}
Consider an $i\in[1,k'-1]$. By
\eqref{eq:hardness-order-coordinates}, the only nonidentity
contributions in the $i$-th coordinate of $H = D^{k'-1}$ from \eqref{eq:hardness-group} come from the elements
 $c_i$ and $c_{i+1}$. Their product is $ab$ if $c_i$ appears in the product before $c_{i+1}$ and $ba$ otherwise. Since the corresponding 
 coordinate of the target $h$ is $ab\neq ba$, $v_{p,i_p}$ must occur before 
 $v_{p+1,i_{p+1}}$ for all $p \in [1,k-1]$, $v_{k,i_k}$ must occur before $w_{1,j_1}$ and 
 $w_{t,j_t}$ must occur before $w_{t+1,j_{t+1}}$ for all $t \in [1,\ell-1]$.
 
 In the $t$-th coordinate ($t\in[1,\ell]$) from the factor $(C_2^\mu)^{\ell}$ in \eqref{eq:hardness-group}, the only non-zero contributions come
from $v_{k_t,i_{k_t}}$ and $w_{t,j_t}$. Since the target $h$ has $\overline{0}$ in this coordinate,
we obtain
\begin{equation}
 \operatorname{code}(i_{k_t})+\operatorname{code}(j_t)=0,
 \text{ i.e. }j_t=i_{k_t}
 \label{eq:hardness-consistency}
\end{equation}
since the mapping code is injective and every element of $C_2^\mu$ is its own inverse.

Fix a $u\in U$ and let $x_p=\chi_{u,i_p}\in\{0,1\}$ for $p\in[1,k]$.
The generators \eqref{eq:vpi} act trivially in the $u$-coordinate of 
$(D_q)^U$ from \eqref{eq:hardness-group}. By the
order of the generators in \eqref{eq:order-of-elements} and \eqref{eq:hardness-consistency}, their product in this coordinate
is 
\[ \prod_{t=1}^{\ell}g_{t,\chi_{u,j_t}} = \prod_{t=1}^{\ell}g_{t,\chi_{u,i_{k_t}}} = \prod_{t=1}^{\ell}g_{t,x_{k_t}}.
\] The target $h$ has the identity in this coordinate, so \eqref{eq:dihedral-test} implies
that $\chi_{u,i_p} = x_p=1$ for at least one $p\in[1,k]$. Thus $u$ is
contained in one of the sets $S_{i_1},\ldots,S_{i_k}$. Since this holds for every
$u\in U$, equation~\eqref{eq:hardness-cover} follows.

For the converse direction, suppose that $i_1,\ldots,i_k\in[1,m]$ satisfy
\eqref{eq:hardness-cover} and consider the product
\begin{equation*} \label{eq:product-converse}
\omega = v_{1,i_1} v_{2,i_2} \cdots v_{k,i_k} w_{1,i_{k_1}} w_{2,i_{k_2}} \cdots w_{\ell,i_{k_\ell}},
\end{equation*}
which has exactly $k+\ell = k'$ factors. The projection of $\omega$ to the factor $C_2^{k'}$ from \eqref{eq:hardness-group}
is $\overline{1}$.
 Moreover, in each coordinate  $t \in [1,\ell]$ of the factor $(C_2^\mu)^\ell$ the product 
 yields $\operatorname{code}(i_{k_t})+\operatorname{code}(i_{k_t})=0$, where
 the first $\operatorname{code}(i_{k_t})$ comes from $v_{k_t, i_{k_t}}$
 and the second $\operatorname{code}(i_{k_t})$ comes from  $w_{t,i_{k_t}}$.
 The projection of $\omega$ to $H$ is
$c_1\cdots c_{k'}=\gamma$. Finally, for each $u\in U$, at least one of the bits
$\chi_{u,i_p}$ for $p \in [1,k]$ is $1$, so \cref{lem:dihedral-test} gives the
identity in the  $u$-coordinate of $D_q^U$. This shows that $\omega$ is the target $h$.

It remains to argue that the reduction is fpt. Let $n_0 = |U| \cdot m + k$ be the size of the parameterized 
\textsc{Set Cover} instance.
The output parameter $k'=k+3\cdot2^k-2$ only depends on the input parameter $k$.
Since $\mu=\lceil\log_2m\rceil$ and $\ell = 3\cdot2^k-2$,
the permutation degree $n$ of $G$ (see \eqref{eq:hardness-degree}) is bounded by 
$2^k \cdot \poly(n_0)$. The number of constructed permutations $v_{p,i}$ and $w_{t,j}$
is also $k' m = 2^k \cdot  \poly(n_0)$.
This shows that our reduction is indeed fpt.
\end{proof}
By \cref{pre:lem:replication,pre:lem:padding}, the problems
$\pF[\mathsf{NilPG}]$, $\pSSS_{\leq}[\mathsf{NilPG}]$, and
$\pSSS[\mathsf{NilPG}]$ are $\W[2]$-hard as well, already for
permutation $2$-groups, since these reductions preserve the group.

\begin{remark}
The group $D_q$ for $q = 2^{k+1}$ is $(k+1)$-step nilpotent: its lower central
series satisfies $\gamma_j(D_q)=\langle r^{2^{j-1}}\rangle$ for
$j\in[2,k+2]$. The group $D$ has nilpotency class $2$, so the
output group $G$ also has nilpotency class $k+1$.
Thus the nilpotency class in this reduction grows
with the source parameter $k$. \cref{thm:nilpotent-hardness}
does not assert $\W[2]$-hardness for any fixed nilpotency class.
\end{remark}

\section{Open problems}

We conclude with some open problems. The most pressing open problems are:
\begin{itemize}
\item the exact complexity of the parameterized subset sum problem over the integers: We lowered 
the best known unconditional upper bound from $\W[3]$ \cite{BussIslam2007} to $\W[2]$. The best known lower bound is $\W[1]$ from \cite{DowneyFellows1995II}.
The conditional upper bound of $\W[1]$ from \cite{AbboudLewiWilliams2014} suggests that the problem is $\W[1]$-complete.
\item the exact complexity of the parameterized factorization problem $\pF[\Sym(*)]$ for symmetric groups: We improved
the previously known lower bound of $\W[1]$ from \cite{CaiChenDowneyFellows1997} to $\W[2]$. The best known upper bound remains $\W^{\func}[2]$ from \cite{LohreyRosowski2025}. 
\end{itemize}

\section{Disclosure of AI use}

GPT-6 Astra produced first versions of the proofs of the four main results \cref{thm:integers}, \cref{thm:solvable}, \cref{thm:nilpotent} and \cref{thm:nilpotent-hardness}. 
The author rewrote the proofs and thereby simplified some details. Apart from this, GPT-6 Astra was used
for literature search. The author takes full responsibility for the final content.


\end{document}